\documentclass[review]{elsarticle}

\usepackage{lineno}
\usepackage{times}
\usepackage{calc}
\usepackage{graphics}
\usepackage{epsfig}
\usepackage{amssymb}
\usepackage{amsmath}
\usepackage{flafter}
\usepackage{float}
\usepackage{latexsym}
\usepackage{psfrag}
\usepackage{dsfont}
\usepackage{rotating}
\usepackage{longtable}
\usepackage{subfigure}
\usepackage{color}
\modulolinenumbers[5]

\journal{Journal of \LaTeX\ Templates}
\renewcommand{\vec}[1]{\mathbf{#1}}

\begin{document}

\begin{frontmatter}

\title{On the use of kinetic energy preserving DG-schemes for large eddy simulation}

\author{David Flad\fnref{myfootnote}}
\address{Institute for Aerodynamics and Gas Dynamics, University of Stuttgart, Pfaffenwaldring 21, 70569 Stuttgart, Germany}
\fntext[myfootnote]{Corresponding Author}

\author{Gregor Gassner}
\address{Mathematical Institute, University of Cologne, Weyertal 86-90, 50923 Cologne, Germany}




\begin{abstract}
Recently, element based high order methods such as Discontinuous Galerkin (DG) methods and the closely related flux reconstruction (FR) schemes have become 
popular for compressible large eddy simulation (LES). Element based 
high order methods with Riemann solver based interface numerical flux functions 
offer an interesting dispersion dissipation behaviour for multi-scale problems: dispersion errors are very low for a broad range of scales, while dissipation 
errors are very low for well resolved scales and are very high for scales close to the Nyquist cutoff. In some sense, the inherent numerical dissipation caused 
by the interface Riemann solver acts as a filter of high frequency solution components. This observation motivates the trend that element based high order methods with 
Riemann solvers are used without an explicit LES model added. Only the high frequency type inherent dissipation caused by the Riemann solver at the element interfaces is used to account for the 
missing sub-grid scale dissipation. Due to under-resolution of vortical dominated structures typical for LES type setups, element based high order methods suffer from stability 
issues caused by aliasing errors of the non-linear flux terms. A very common strategy to fight these aliasing issues (and instabilities) is so-called polynomial de-aliasing, where interpolation is exchanged with projection based on an increased number of quadrature points. In this paper, we start with this common no-model or implicit LES (iLES) DG approach with polynomial de-aliasing and Riemann solver dissipation and 
review its capabilities and limitations. We find that the strategy gives excellent results, but only when the resolution is such, that about 40\% of the dissipation is resolved. For more realistic, coarser resolutions used in classical LES e.g. of industrial applications, the iLES DG strategy becomes quite in-accurate. We show that there is no obvious fix to this strategy, as adding for instance a sub-grid-scale models on top doesn't change much or in worst case decreases the fidelity even more. Finally, the core of this work is a novel LES strategy based on split form DG methods that are kinetic energy preserving. Such discretisations offer 
excellent stability with full control over the amount and shape of the added artificial dissipation. This premise is the main idea of the work and we will assess the LES capabilities of the novel split form DG approach. We will demonstrate that the novel DG LES strategy offers similar accuracy as the iLES methodology for well resolved cases, but strongly increases fidelity in case of more realistic coarse resolutions. 
\end{abstract}


\end{frontmatter}

\section{Introduction}
\label{sec:intro}
Using high order Discontinuous Galerkin (DG) methods, and closely related methods such as flux reconstruction (FR) schemes, implicit large eddy simulation 
(iLES) approaches were recently applied to successfully simulate flows at moderate Reynolds numbers e.g. 
\cite{uranga2011,wiart2012,les1,les2,laf,wiart2015,moura2016}. The approach used in these references does not add an explicit sub-grid scale model, 
such as e.g. eddy viscosity, but instead relies on the dissipation behaviour of the schemes due to Riemann solvers at element interfaces. 
High order discretisations, i.e. methods with polynomial degree $N\geq 3$, are usually used to obtain good dispersion behaviour which is necessary to accurately 
capture the interaction of different spatial scales. Within the references above, it is shown that for very high order polynomial degrees, the iLES DG discretisations are prone to 
aliasing issues due to the non-linearity of the flux functions. These aliasing issues can even cause instabilities and cause the simulation to crash. Within the references above, 
stabilisation in case of very high order polynomial degrees is done by applying so-called polynomial de-aliasing. There are several ways of 
implementing and interpreting polynomial de-aliasing. The goal however is always the same: account for the non-linearity 
of the flux function by either using a polynomial projection or by increasing the number of quadrature points for the approximation of the integrals. 
In  \cite{Kirby2003}, Kirby and Karniadakis demonstrated that cutting off $1/3$ of the highest modes in every time step for the incompressible Navier-Stokes 
equations (quadratic flux function) gives a stable approximation. They extended this for compressible flow, where they interpreted the non-linearity of the 
compressible Navier-Stokes fluxes as cubic and found that cutting off $1/2$ of the highest modes gives stability in all cases they tested.  The non-linearity in 
the compressible case is however rational and thus there are also cases where this strategy fails, see e.g. \cite{moura2016,gassner2016split}. The issue of 
aliasing is mitigated somewhat by the use of Riemann solvers at element interfaces. For element based high order methods, the Riemann solver at the 
element interfaces provides numerical dissipation that acts in a high frequency filter-like way: it is almost zero for well resolved scales but gets very high for scales close 
to the Nyquist limit. This numerical dissipation behaviour paired with the good dispersion properties is the motivation to use element based high order methods 
such as the DG scheme with interface Riemann solvers for under-resolved turbulence simulations, often termed implicit LES (iLES), as no explicit sub-grid scale 
dissipation model is added. 

\begin{figure}[!htpb]
  \subfigure[]{%
    \label{fig:highresOI}
    \includegraphics[width=0.5\textwidth, height=4.5cm, trim=20 10 60 20 ,clip]{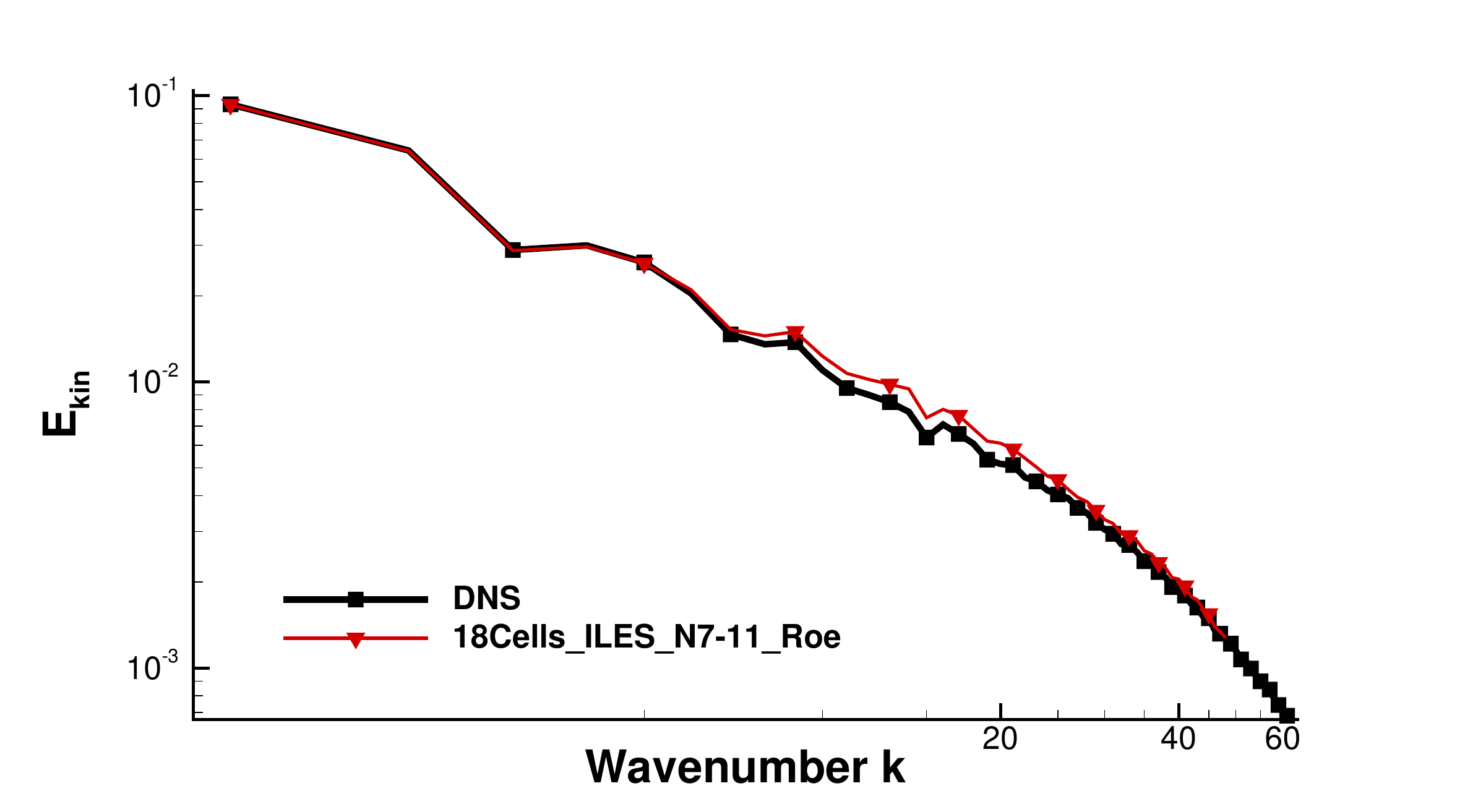}}
  \subfigure[]{%
    \label{fig:lowresOI}
    \includegraphics[width=0.5\textwidth, height=4.5cm, trim=20 10 60 20 ,clip]{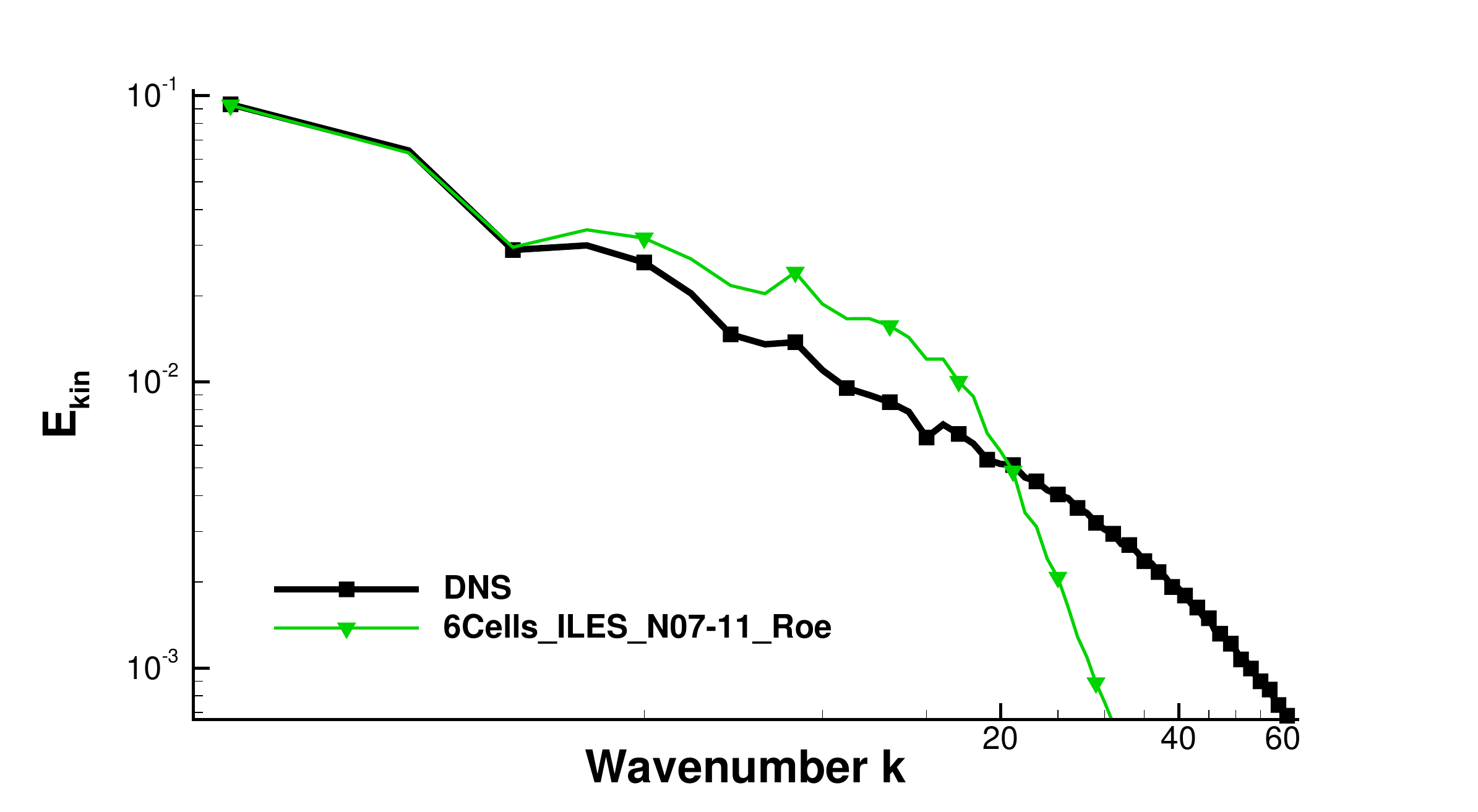}}\\
  \caption{Decaying homogeneous isotropic turbulence DNS/LES kinetic energy spectra: a) iLES $144^3$ DOF b) iLES $48^3$ DOF. Both approximations use Roe's 
approximate Riemann solver at the element interfaces.
}
\end{figure}

In the following, we will use the iLES DG approach to simulate the decay of homogeneous isotropic turbulence (DHIT). The test case is a homogeneous isotropic 
turbulence freely decaying for about one large eddy turn over time ($T=\bar{v}(t_{start})/L_{int}\approx1.3$). The Reynolds number based on the Taylor micro scale decays from 
$Re_{\lambda}\approx162$ to $Re_{\lambda}\approx97$. The LES initial state is obtained from a filtered DNS field. The 
DNS is simulated as in \cite{yamazaki2002effects}, with a pseudo-spectral code on $512^3$ DOF using the $2/3$ rule for de-aliasing. Similar test cases for LES 
were used for example in \cite{VMShughes2000,hickel2006,laf}. Figure \ref{fig:highresOI} shows the resulting kinetic energy (KE) spectra for a relatively high 
resolved LES setting, $18$ cells per direction and polynomial degree $N=7$, i.e. $144$ DOF per direction in comparison with the DNS result.  For the Riemann 
solver, we choose Roe's approximate numerical flux function. We note that in this case polynomial de-aliasing is applied according to the $3/2$ rule (total of $12^3$ quadrature points per element), as the flow behaves 
almost incompressible, see e.g. \cite{tcfd2012}.

 With this particular iLES DG setup about $62.5\%$ of the total dissipation at $t_0$ can be resolved if we consider 
the theoretical Nyquist wavenumber $k_{Ny}=144/2=72$ for this discretisation with $144^3$ DOF. Assuming a more realistic (but still optimistic) 
approximation capability of the piece-wise polynomial ansatz space of $3$ points per wavenumber (PPW), then only about about $43.2\%$ of the KE dissipation is 
resolved. This setup is still somewhat well behaved: although under-resolved, the cut off wavenumber $k_c=144/3=48$ falls within the beginning of the 
dissipation range, which is the reason that the results obtained with the DG iLES approach are in excellent agreement with the reference direct numerical 
simulation (DNS).

However, in the available literature, simulations of the DHIT test case usually aim at considerable lower resolutions. In the spectral community $32^3$ points, e.g. \cite{park2004discretization}, are 
commonly used for the assessments. Accounting for the difference between approximations with Fourier basis ($2$ PPW resolution limit) and polynomial basis (about $3$ PPW 
resolution limit), we increase the DG resolution to $48^3$ DOF, in particular we choose for the next setup $6^3$ elements with a polynomial degree of $N=7$. Using again the 
optimistic estimate of $3$ PPW the resulting cut-off wavenumber of this setup is $k_c=48/3=16$, which is comparable to the cut-off wavenumber of a spectral 
simulation with $32^3$ DOF, i.e. $k_c=32/2=16$. If we estimate the resolved total dissipation for this coarse setup, we only get about $0.9\%$ for the cut-off wavenumber $k_c=16$, which shows 
that this test case is severely under-resolved in comparison to the setup with $144^3$ DOF. Figure \ref{fig:lowresOI} shows that for such a resolution the DG iLES 
fails to match the reference DNS result. We observe a pile up of energy for mid range scales followed by a strong drop off of energy. Such a behaviour is indeed 
typical for upwind based iLES discretisations with very coarse resolutions, such as the present DG iLES configuration with Roe's Riemann solver. A discussion on this typical behaviour of upwind schemes and 
a good overall review can be found in Drikakis \cite{drikakis2009}.

\begin{figure}[!htpb]
    \centering
    \includegraphics[width=0.5\textwidth, height=4.5cm, trim=20 10 60 20 ,clip]{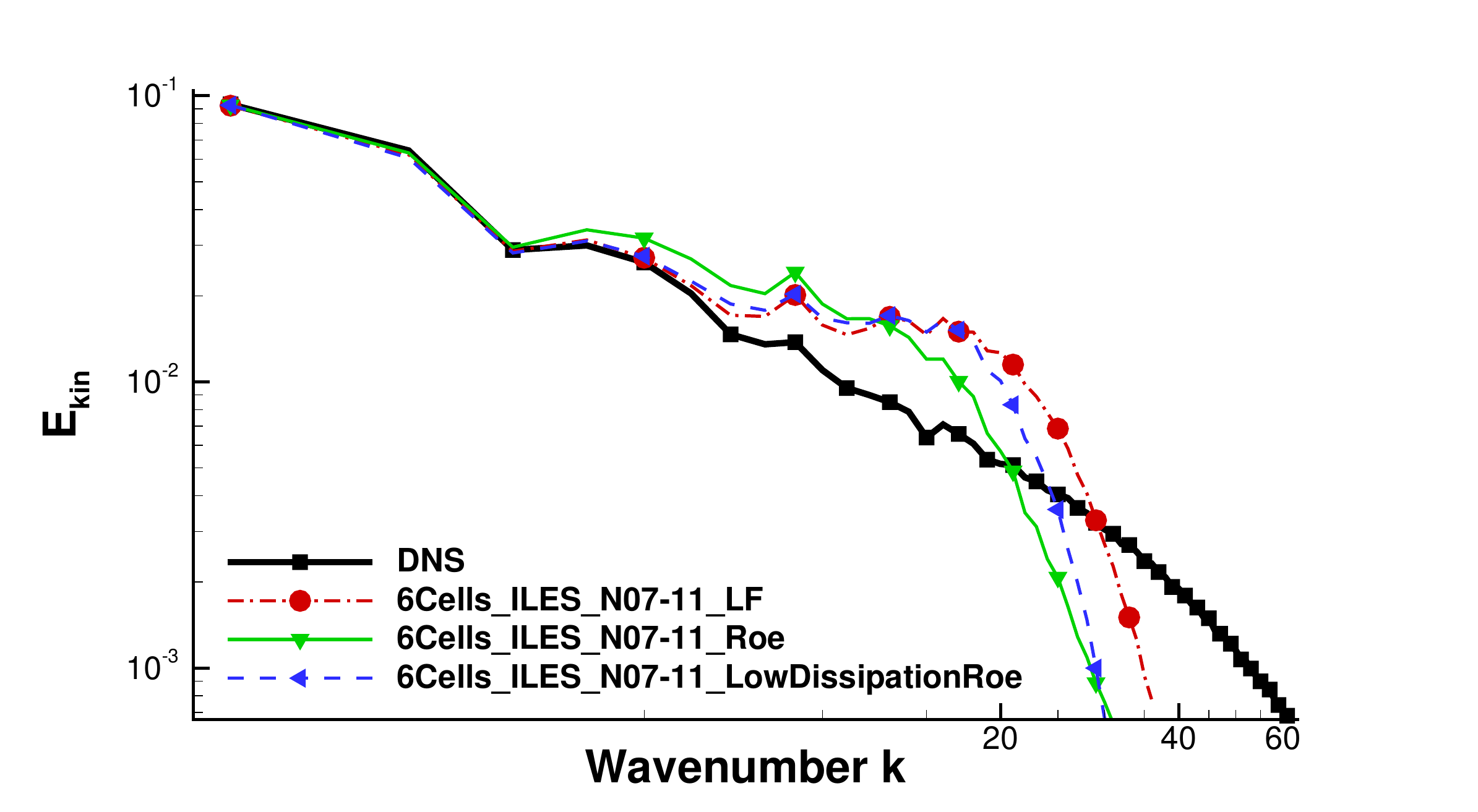}
  \caption{Decaying homogeneous isotropic turbulence DNS/LES kinetic energy spectra: iLES with different Riemann solvers. Red curve is the result with 
Lax-Friedrichs and green with Roe's flux.}
  \label{fig:iLESRiemann}
\end{figure}

As discussed, for the DG iLES approach the main dissipation is added by the Riemann solver at the element interfaces. Consequently, the amount and the 
'shape' of the numerical dissipation depends on the particular choice of the Riemann solver. Figure \ref{fig:iLESRiemann} shows results when comparing Roe's 
approximate Riemann solver and the Lax-Friedrichs numerical flux. While there is a clear effect and influence by the choice of the interface numerical flux 
function, neither choice greatly improves the iLES result in comparison to the DNS reference, i.e. both suffer from the typical upwind behaviour mentioned above. 

\begin{figure}[!htpb]
    \centering
    \includegraphics[width=0.5\textwidth, height=4.5cm, trim=20 10 60 20 ,clip]{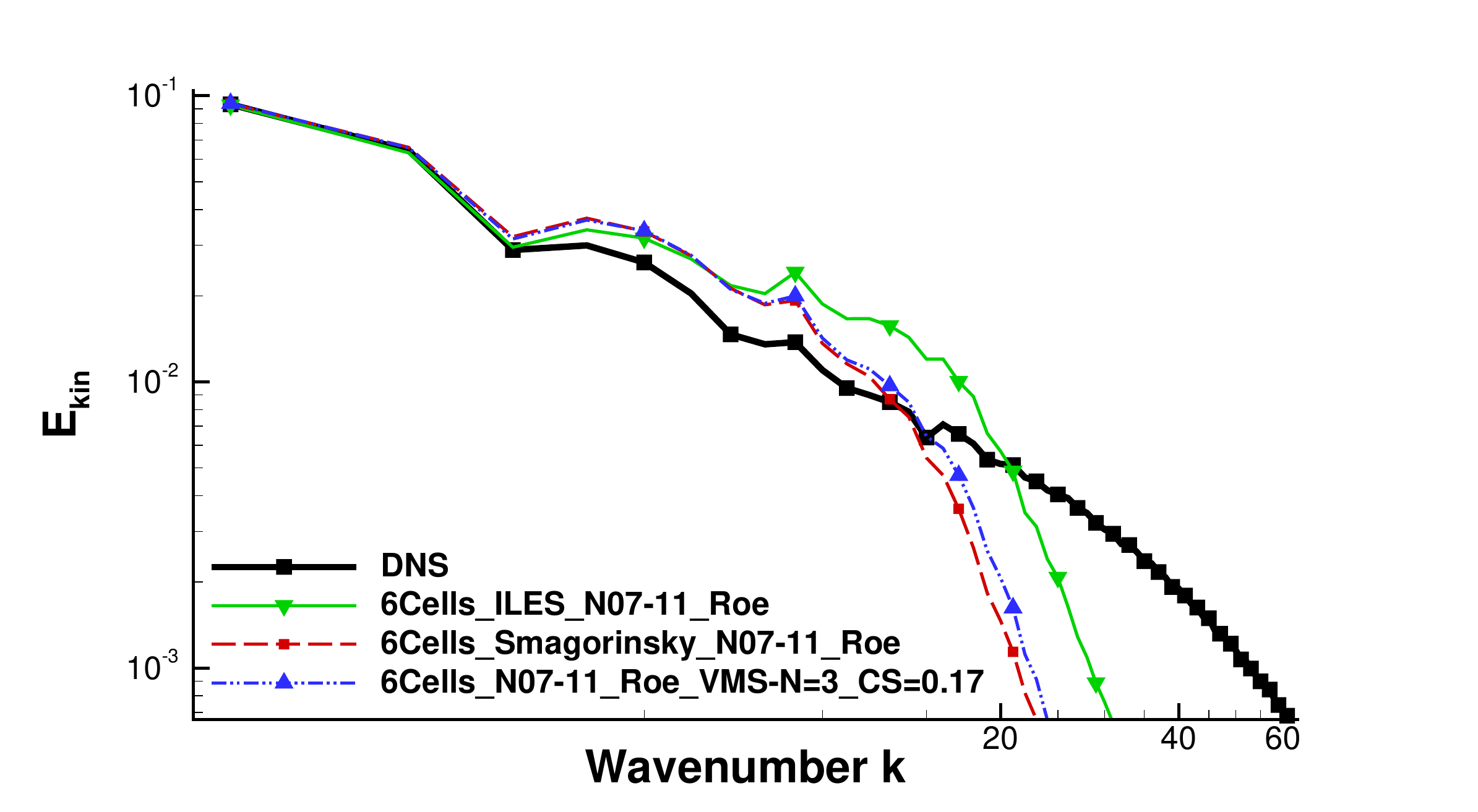}
  \caption{Decaying homogeneous isotropic turbulence DNS/LES kinetic energy spectra: polynomial de-aliasing, Roe's Riemann solver and two different explicit 
sub grid scale eddy viscosity models}
  \label{fig:iLESRiemann_expl}
\end{figure}

It is worth noting that the second 'parameter' in the DG iLES approximation, the polynomial de-aliasing, cannot be omitted as the DG iLES without polynomial 
de-aliasing becomes unstable for this particular choice of polynomial degree (it might be stable for lower polynomial degrees, which results in higher artificial upwind dissipation). 
In conclusion, it seems that the DG iLES approach is very 'rigid', rigid in the sense that there is not a lot of options to 
influence the dissipation behaviour of the scheme as polynomial de-aliasing is necessary and the only available parameter 'Which Riemann solver do we chose?' doesn't seem 
to provide means for substantially increasing the fidelity of the DG iLES. Adding explicit SGS dissipation to the described (iLES) schemes is difficult as 
usually the artificial dissipation interacts with the sub-grid scale model viscosity. It is then difficult to separate the 
effects, making turbulence modelling based on physical rationals very adventurous. This could be in fact the 
explanation, why the no-model or iLES DG approach is current state of the art in the DG community, as not really a clearly better strategy is known yet. Exemplary we show the effect of 
adding an explicit SGS model in Figure \ref{fig:iLESRiemann_expl} for the standard Smagorinsky model, see \ref{App:smago} for details, and its high fidelity variant, 
the small-small variational multi scale model (VMS), see \ref{App:VMS} for details.
We find again that indeed there is no significant improvement of the result, independent of the SGS model constants used. This strongly hints towards a limiting applicability of upwind-type based iLES DG with de-aliasing for cases with relatively high resolution: Our numerical experiments show that roughly $40\%$ of the total dissipation should be resolved, which poses a natural plateau for 
this approach.Up to now, it seems that no remedy for this applicability plateau is available and known in literature. That is the reason why in this work we do not propose to modify or add models to 
the DG iLES to fix these issues, but propose a completely new approach instead! 

We change the discretisation considerably and focus on an alternative approach for stabilising under-resolved DG discretisations, namely the strategy of using 
split forms for the non-linear terms \cite{gassner2016split}, e.g. skew-symmetric splittings of the advective terms. The class 
of split form DG allows to inbuilt the de-aliasing in the discretisation of the volume terms. It is furthermore possible to construct a high order DG 
discretisation that is kinetic energy preserving, if an appropriate split form is chosen. Split forms are well known in the finite difference community 
\cite{park2004discretization,pirozzoli2011numerical,nagarajan2003robust,mittal1997suitability,kravchenko1997effect,desjardins2008high} and a multitude of 
variants exist. A good overview is given in Pirozzoli \cite{pirozzoli2011numerical}. Priozzoli also presents a special variant of the splitting based on a 
three-way splitting of the non-linear flux terms, resulting in a kinetic energy preserving discretisation that is close to the scheme presented by Kennedy and 
Gruber \cite{kennedy2008}. Based on these split forms, it is possible to construct dissipation free DG schemes that still feature remarkable numerical robustness, 
sometimes even higher than the polynomial de-aliasing based approach \cite{gassner2016split}. This feature is exactly what we will exploit in this work: having 
a dissipation free base line scheme that is still robust offers the opportunity to add precisely the shape and amount of dissipation motivated by modelling sub-grid scale physics instead 
of e.g. Riemann solver dissipation that is constructed to capture shock physics.

Summarising, the goal in this work is to show how to construct dissipation free split form DG schemes with explicit sub-grid turbulence models and assess the performance of this novel strategy. We demonstrate 
that this strategy offers high potential for practical LES applications were only a small fraction of the total dissipation is resolved. We further demonstrate that it is possible to overcome the applicability plateau of the DG iLES approach and that with a dissipation free base line scheme we are able to add dissipation motivated by turbulence-physics. 

The remainder of the article is organised as follows: In the next section, a small introduction to split form DG schemes is given. The next two chapters focus on numerical results for the DHIT with the conclusions and outlook given in the last section.

\section{Numerical Methods}
\label{sec:numerical_method}

\subsection{Discontinuous Galerkin Spectral Element Method}
\label{sec:dgsem}

We consider the compressible Navier-Stokes equations (NSE) expressed
in conservation form 
\begin{equation}\label{eq:NSE}
\vec{U}_t + \vec{\nabla}_{x}\cdot\vec{F}(\vec{U},\vec{\nabla}_x \vec{U}) = 0,
\end{equation}
where $\vec{U}$ denotes the vector of conserved quantities $\vec{U}=(\rho, \rho v_1, \rho
v_2, \rho v_3, \rho e)^T$, the subscript $t$ the time derivative and $\vec{\nabla}_x$ the gradient operator in physical space. The flux
is the difference of advection and viscous fluxes, $\vec{F}=\vec{F}^{a}(\vec{U})
- \vec{F}^{v}(\vec{U},\vec{\nabla}_{x}\vec{U})$, with the entries 
\begin{equation}
\vec{F}^{a}_{l}(\vec{U})=\begin{pmatrix}\rho\, v_l\\ \rho\, v_1 v_l +\delta_{1l}\,p\\
\rho\, v_2 v_l +\delta_{2l}\,p\\ \rho\, v_3 v_l +\delta_{3l}\,p  \\ \rho\,ev_l + p\,v_l\end{pmatrix},\; 
\vec{F}^{v}_{l}(\vec{U},\vec{\nabla}_{x}\vec{U})=\begin{pmatrix}0\\
\tau_{1l}\\
\tau_{2l}\\\tau_{3l}\\ \tau_{lj}v_j - q_l\end{pmatrix},
\end{equation}
where $l=1,2,3$, denoting the Cartesian directions of the flux $\vec{F}_1,\vec{F}_2,\vec{F}_3$. We follow the usual nomenclature for $\rho, 
(v_1,v_2,v_3)^T, p,
e$ denoting the density, velocity vector, pressure and specific total energy, respectively.
The perfect gas law
\begin{equation}
p=\rho R T=(\kappa-1)\rho(e-\frac 1 2
\vec{v}\cdot\vec{v}),\quad
e=\frac{1}{2}\vec{v}\cdot\vec{v} + c_v T
\end{equation}
is used to close the system of equations.

In this work we use a special DG variant, namely the discontinuous Galerkin spectral element collocation method (DGSEM) with Legendre Gauss-Lobatto (LGL) nodes. 
The LGL nodes are essential, as this choice guarantees the so-called summation-by-parts (SBP) property of the resulting DGSEM operator 
\cite{gassner_skew_burgers}. Up to now, split form DG is only available for this specific variant, as the SBP property is fundamental. The computational domain 
is subdivided into 
non-overlapping hexahedral elements which are transformed to reference space $\vec{\xi}$ via a transfinite mapping. Within the reference element, a tensor-product polynomial approximation is constructed: we use the tensor-product of the 1D LGL nodes and the associated tensor-product of 1D Lagrange polynomials, which gives $(N+1)^3$ DOF per element per unknown quantity for a given polynomial degree of $N$:
\begin{equation}
 \label{eq:ansatz}
 \vec{U}(\vec{\xi},t)\approx\sum\limits_{i,j,k=0}^N \vec{\hat{U}}_{ijk}(t)\psi_{ijk}(\vec{\xi})\,,\qquad 
\psi_{ijk}(\vec{\xi})=\ell_i(\xi^1)\ell_j(\xi^2)\ell_k(\xi^3)\,,
\end{equation}
with the 1D Lagrange basis functions defined as
\begin{equation}
 \ell_j(\vec{\xi}) = \prod_{i=0,i\neq j}^N\frac{\xi-\xi_i}{\xi_j-\xi_i}, \quad j=0,...,N,
\end{equation}
where $\{\xi_i\}_{i=0}^N$ are the LGL nodes. From this 1D Lagrange basis, we can directly compute the associated derivative matrix 
\begin{equation}
D_{ij} = \ell_j'(\vec{\xi}_i),\quad i,j=0,...,N,
\end{equation}
store the LGL weights $\{\omega_i\}_{i=0}^N$ in the discrete mass-lumped diagonal mass matrix 
\begin{equation}
M = diag([\omega_0,...,\omega_N]),
\end{equation}
and define the boundary matrix
\begin{equation}
B= diag([-1,0,...,0,1]).
\end{equation}
These three 1D operator matrices are connected via the SBP property
\begin{equation}
(M\,D) + (M\,D)^T = B,
\end{equation}
which guarantees consistency of discrete integration-by-parts based on the derivative matrix and the discrete inner product with LGL nodes. 

For simplicity, we assume a Cartesian grid as used in our numerical tests for the DHIT test case. Details and the extension to general curvilinear elements are 
available in \cite{gassner2016split}. Furthermore, for the discretization of the second order viscous terms we resort to a standard method available in 
literature, namely, we use the approach introduced by Bassi and Rebay (BR) \cite{Bassi&Rebay:1997:B&F97,bassi2} and analysed in Gassner et al. \cite{br1isstable}, where the gradient of the solution is 
introduced as an additional new unknown. With this, the second order problem is re-written into a larger system of first order partial differential equations. 
For the system of auxiliary gradients a standard DGSEM based on LGL is used, where the numerical flux at the element interfaces is chosen as the 
arithmetic mean. This DG gradient is then used to compute the viscous fluxes in the discrete Navier-Stokes equations. There are two different variants which can 
be used in the BR framework, namely BR1 (arithmetic mean of the viscous fluxes at the element interface) or BR2 (an additional artificial dissipation term 
depending on the interface jump of the solution). The BR2 scheme provides a mechanism to introduce additional numerical dissipation via the so-called penalty 
constant $\eta_{BR}\geq 1$: the larger this value, the higher the additional artificial interface dissipation.

The discretization of the advective Euler fluxes needs more care. Assuming Cartesian grids with size $\Delta x, \Delta y, \Delta z$ the transformation of the 
flux divergence in reference space is straight forward. For simplicity, we only focus on the term $\frac{1}{\Delta x} \vec{F}^a_1(\vec{U})_\xi$, i.e. the 
Cartesian flux in 
$\xi$-direction after transformation. The standard strong form LGL DGSEM for this term at an LGL node $(i,j,k)$ reads as
\begin{equation}
\frac{1}{\Delta x} \vec{F}^a_1(\vec{U})_\xi\big|_{ijk} \approx  \frac{1}{M_{ii}}\left(\delta_{iN}\left[\vec{F}_1^{a,*} - \vec{F}_1^a\right]_{Njk} - 
\delta_{{i0}}\left[\vec{F}_1^{a,*} - \vec{F}^a_1\right]_{0jk}\right)+\sum_{m=0}^N D_{i\boldsymbol{m}}(\vec{F}^a_1)_{\boldsymbol{m}jk},
\end{equation}
where $\vec{F}_1^{a,*}$ is the numerical interface flux function, typically an approximate Riemann solver, and $(\vec{F}^a_1)_{mjk} 
=\vec{F}^a_1(\boldsymbol{\hat{U}}_{ijk})$ is the non-linear Euler flux evaluated at the nodal solution value, i.e. the interpolation of the Euler flux at the 
LGL grid. The Kronecker delta $\delta_{ij}$ is equal to one when $i=j$ and zero otherwise and acts as a switch such that the boundary terms only affect the LGL 
boundary nodes. 

As explained above, due to the non-linearity of the Euler fluxes, standard high order DG implementations such as the LGL DGSEM are prone to aliasing issues. An approach discussed above is to use polynomial de-aliasing, where the non-linear flux functions are projected onto the space of polynomials of degree $N$, instead of using an interpolation only. However, due to the rational non-linear nature of the Euler fluxes with respect to the conservative variables, it is not possible to implement an exact projection. Instead, this projection is constructed with an 'over-integration' of the flux terms, i.e. by using a higher number of quadrature nodes. In our case we use the integer of $Q=3/2 (N+1) -1$ as we assume weak compressibility effects due to the relatively low Mach number. A possible way of implementing this is via modal cut-off filters
\begin{equation}
\frac{1}{\Delta x} \vec{F}^a_1(\vec{U})_\xi\big|_{ijk} \approx  P_{Q\to N}\frac{1}{M_{ii}}\left(\delta_{iQ}\left[\vec{F}_1^{a,*} - \vec{F}_1^a\right]_{Qjk} - 
\delta_{{i0}}\left[\vec{F}_1^{a,*} - \vec{F}^a_1\right]_{0jk}\right)+\sum_{m=0}^Q D_{i\boldsymbol{m}}(\vec{F}^a_1)_{\boldsymbol{m}jk},
\end{equation}
where the modal filter $P_{Q\to N}$ transforms the nodal space into modal space (tensor-product of orthogonal Legendre polynomials) deletes the modes from $N+1$ to $Q$ and transforms back to nodal space, see  e.g. \cite{laf} for details.  

The key of this work is an alternative technique for de-aliasing, namely the split form DG methodology. Split form DG is based on a remarkable property of diagonal norm SBP operators, discovered and introduced by Carpenter and Fisher, e.g. \cite{fisher2013,carpenter2014}. They constructed so-called entropy stable forms of the volume terms by using 
\begin{equation}
\frac{1}{\Delta x} \vec{F}^a_1(\vec{U})_\xi\big|_{ijk} \approx  \frac{1}{M_{ii}}\left(\delta_{iN}\left[\vec{F}_1^{a,*} - \vec{F}_1^a\right]_{Njk} - 
\delta_{{i0}}\left[\vec{F}_1^{a,*} - \vec{F}^a_1\right]_{0jk}\right)+\sum_{m=0}^N 2\, D_{i\boldsymbol{m}}\vec{F}^{a,\#}_1(\vec{U}_{ijk}, 
\vec{U}_{\boldsymbol{m}jk}),
\end{equation}
where they used a two-point entropy conserving numerical volume flux $\vec{F}^{a,\#}_1(\vec{U}_{ijk}, \vec{U}_{\boldsymbol{m}jk})$. In Gassner et al. 
\cite{gassner2016split} it was shown that it is possible to choose every symmetric and consistent two-point numerical volume flux $\vec{F}^{a,\#}_1$ and that 
every choice results in a novel split form DGSEM. In particular, the authors identified specific choices of numerical volume fluxes to exactly reproduce well 
known split forms, such as the one of Pirozzoli \cite{pirozzoli2011numerical}. This particular numerical volume flux is used in the present work and reads as
  \begin{equation}
  \label{eq:PIflux}
      \vec{F}^{a,\#}_1(\vec{U}_{ijk},\vec{U}_{mjk})=\begin{bmatrix}
                            \{\{\rho\}\}\{\{u\}\} \\
                            \{\{\rho\}\} \{\{u\}\}^2 + \{\{p\}\} \\
                            \{\{\rho\}\} \{\{u\}\} \{\{v\}\} \\
                            \{\{\rho\}\} \{\{u\}\} \{\{w\}\} \\
                            \{\{\rho\}\} \{\{u\}\} \{\{h\}\}
                            \end{bmatrix},
   \end{equation}
   with
 \begin{equation}\nonumber
    \{\{\alpha\}\} := \frac{1}{2}\,(\alpha_{ijk}+\alpha_{mjk}).
   \end{equation}
 In \cite{gassner2016split}, it was also shown that the resulting split form DG scheme is kinetic energy preserving in the element and when choosing the element 
interface flux $\vec{F}_1^{a,*}$ equal to \eqref{eq:PIflux}, the resulting split form DGSEM is kinetic energy preserving across the domain. Effectively, this means that the aliasing error in the kinetic energy preservation due to the discretisation of advective terms is eliminated. This split form kinetic energy preserving scheme is exactly our 
baseline scheme with zero artificial kinetic energy dissipation mentioned above. Starting with this 
dissipation free baseline method, it is possible to add Riemann solver type dissipation at the element interfaces by augmenting the central numerical flux 
$\vec{F}_l^{a,\#}$ at the interface with either Lax-Friedrichs type scalar dissipation or Roe type matrix dissipation to recover an unwinding type split form DG scheme. However, this baseline scheme offers many other ways of adding artificial dissipation, as we will discuss below in detail. 

Finally, the semi discrete (split form) DGSEM is integrated in time with an explicit fourth order low storage Runge-Kutta method, \cite{Carpenter-MH:2005kx}. The explicit time step is computed on the fly by a typical advective CFL condition.

\section{Results for decaying isotropic turbulence}
\begin{figure}[!htpb]
    \centering
    \includegraphics[width=0.5\textwidth, height=4.5cm, trim=20 10 60 20 ,clip]{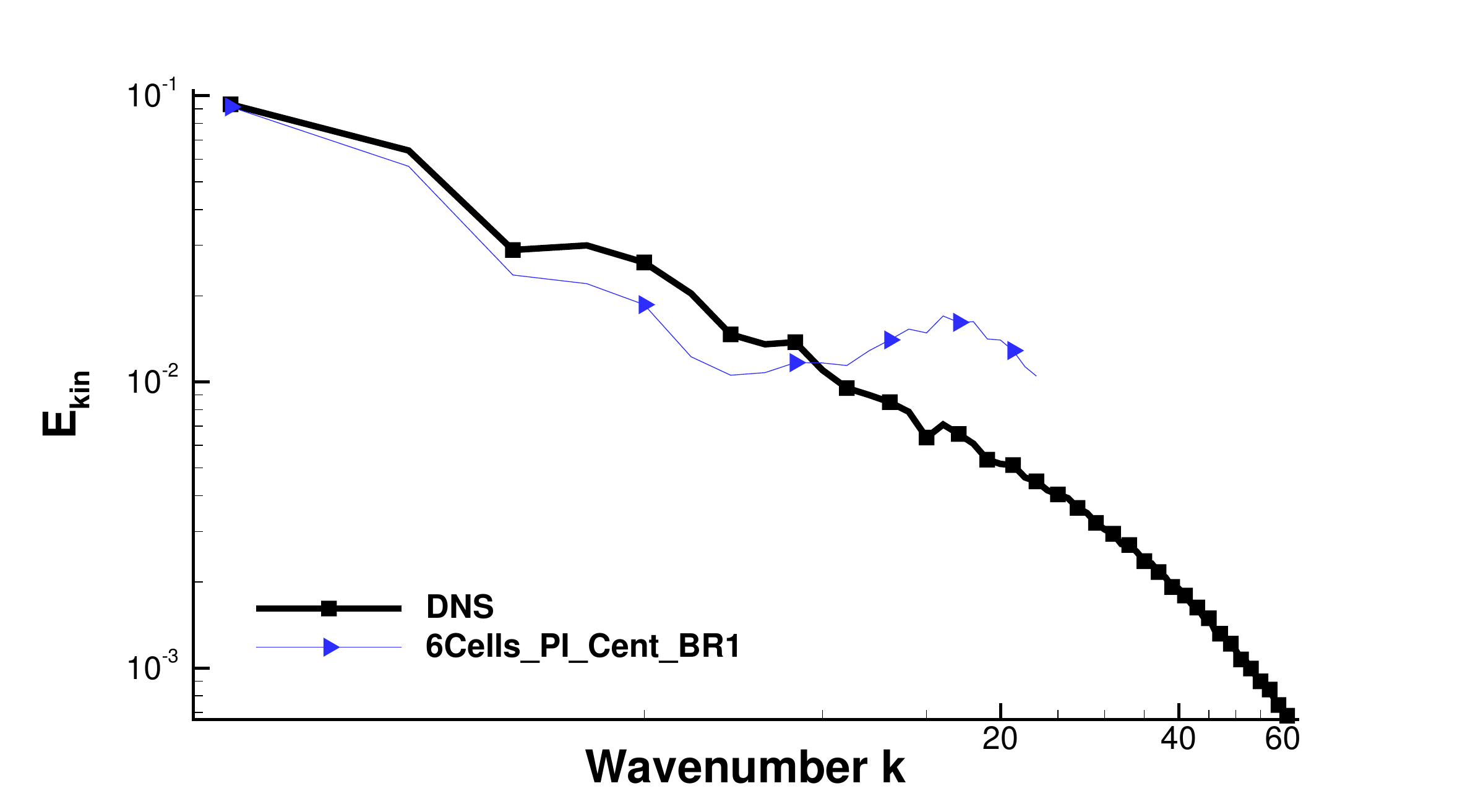}
  \caption{Decaying homogeneous isotropic turbulence DNS/LES kinetic energy spectra: no-model, central (``dissipation free``) numerical fluxes, BR1 scheme. 
Discretization with $6^3$ cells and $N=7$.}
  \label{fig:split_cent}
\end{figure}

Figure \ref{fig:split_cent} shows the no model result for the DHIT test case for a dissipation free scheme employing the flux splitting with the numerical flux of Pirozzoli \eqref{eq:PIflux} discussed above (labeled PI in 
the following figures). 

As seen in this plot the discrete solution describes approximately a $k^2$-slope close to the cut of wavenumber ($\approx 16$), known as equipartition 
spectrum. That is indeed an expected result and a direct consequence of the missing dissipation of the small scales. A result often observed when 
using no model pseudo spectral codes, eg. \cite{park2004discretization}. To the best of our knowledge that is the first time such a result is obtained with a 
DG scheme. It is furthermore remarkable that the dissipation free high order split form DGSEM with $6^3$ elements and $N=7$ is stable, although the 
simulation is severely under-resolved. Besides these two remarkable observations, it is also clear that this scheme does not provide an accurate approximation 
of the reference DNS result. However, as stated above, we will use this dissipation free scheme as a baseline scheme. The goal is to show that this scheme 
offers great potential when adding dissipation for sub-grid scale turbulence modelling. In fact, we will demonstrate several different strategies to include 
sub-grid scale dissipation. 

Lastly, we note that in Figure \ref{fig:split_cent} (and all other figures containing kinetic energy spectra for LES) the spectra of the LES is plotted up to the theoretical Nyquist frequency. The resolution limit of polynomials however is theoretically limited to $\pi$ PPW. For the dissipation free scheme in Figure \ref{fig:split_cent} this limit corresponds well with the start where the kinetic energy 
decreases again ($48$ DOF $\widehat{=}$ $k_{\pi PPW}=15.28$, $144$ DOF $\widehat{=}$ $k_{\pi PPW}=45.84$). In the following sections, for our LES results this 
polynomial cut-off wavenumbers are the natural limits of our approximation. Our goal is therefor to obtain the best possible match with the DNS reference within 
this limit. 

\subsection{Split form DG iLES}
 
Our first strategy is to use an iLES type strategy to account for the missing scales. We show, that the de-aliasing strategy (split form or polynomial de-aliasing) affects the behaviour of the iLES method, however the failure of iLES DG for coarse resolutions is not affected by the choice of de-aliasing strategy. We discuss in the following the deficiencies of the iLES DG approach and offer suggestions that the reason is missing dissipation for the large scales. 

Starting with the baseline scheme that is dissipation free \ref{fig:split_cent}, we first add a classic Roe type matrix dissipation term at the element interfaces, which adds an upwind-type dissipation to the scheme based on linearised characteristics. This additional dissipation implicitly models the missing dissipation of the sub-grid scales, analogous to the standard DG iLES approach with polynomial de-aliasing discussed in the introduction. The core difference now is that we completely changed the de-aliasing mechanism and are thus able to investigate its effect on the iLES DG performance.

\begin{figure}[!htpb]
  \subfigure[]{%
    \label{fig:iLES12split}
    \includegraphics[width=0.5\textwidth, height=4.5cm, trim=20 10 60 20 ,clip]{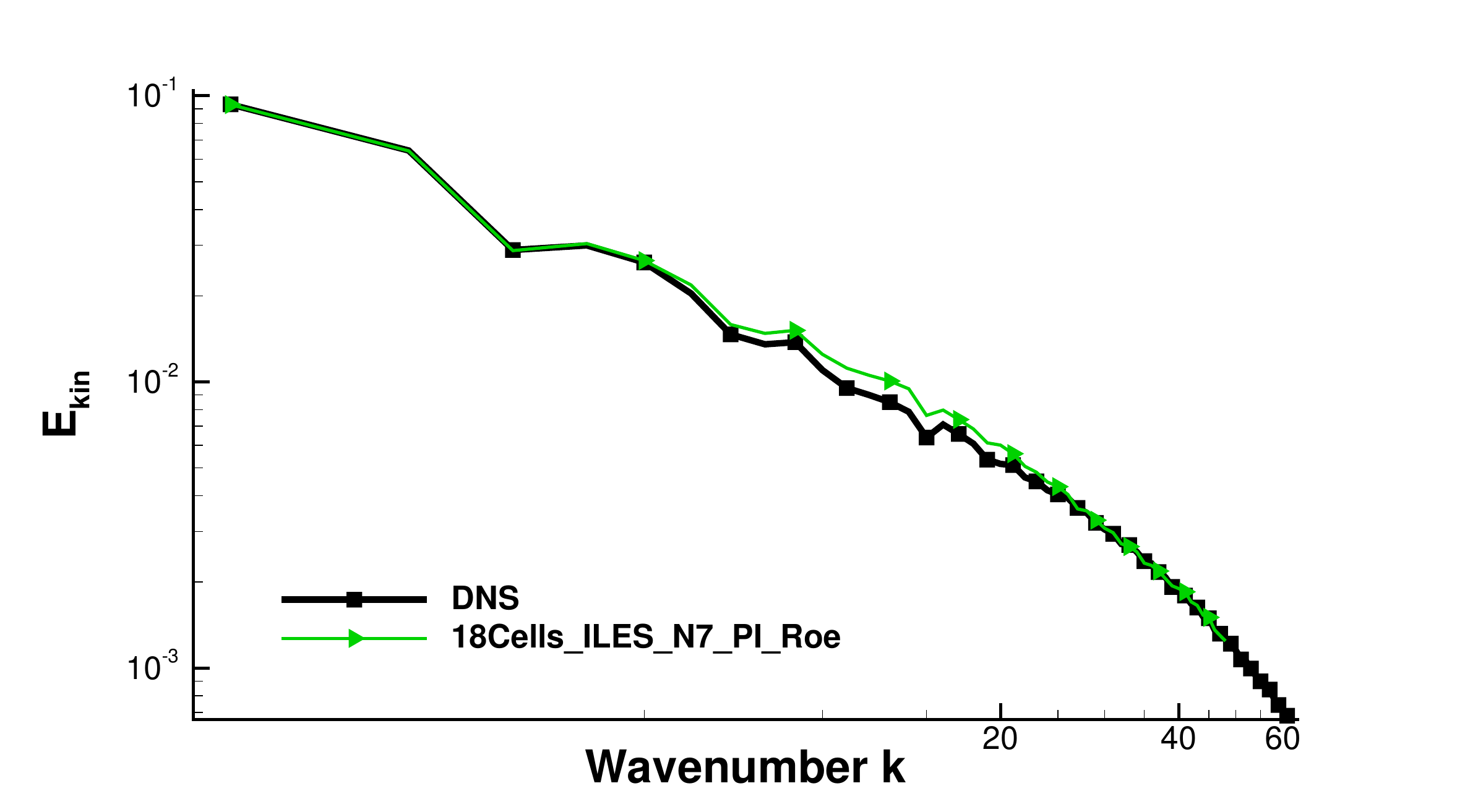}}
  \subfigure[]{%
    \label{fig:iLES6split}
    \includegraphics[width=0.5\textwidth, height=4.5cm, trim=20 10 60 20 ,clip]{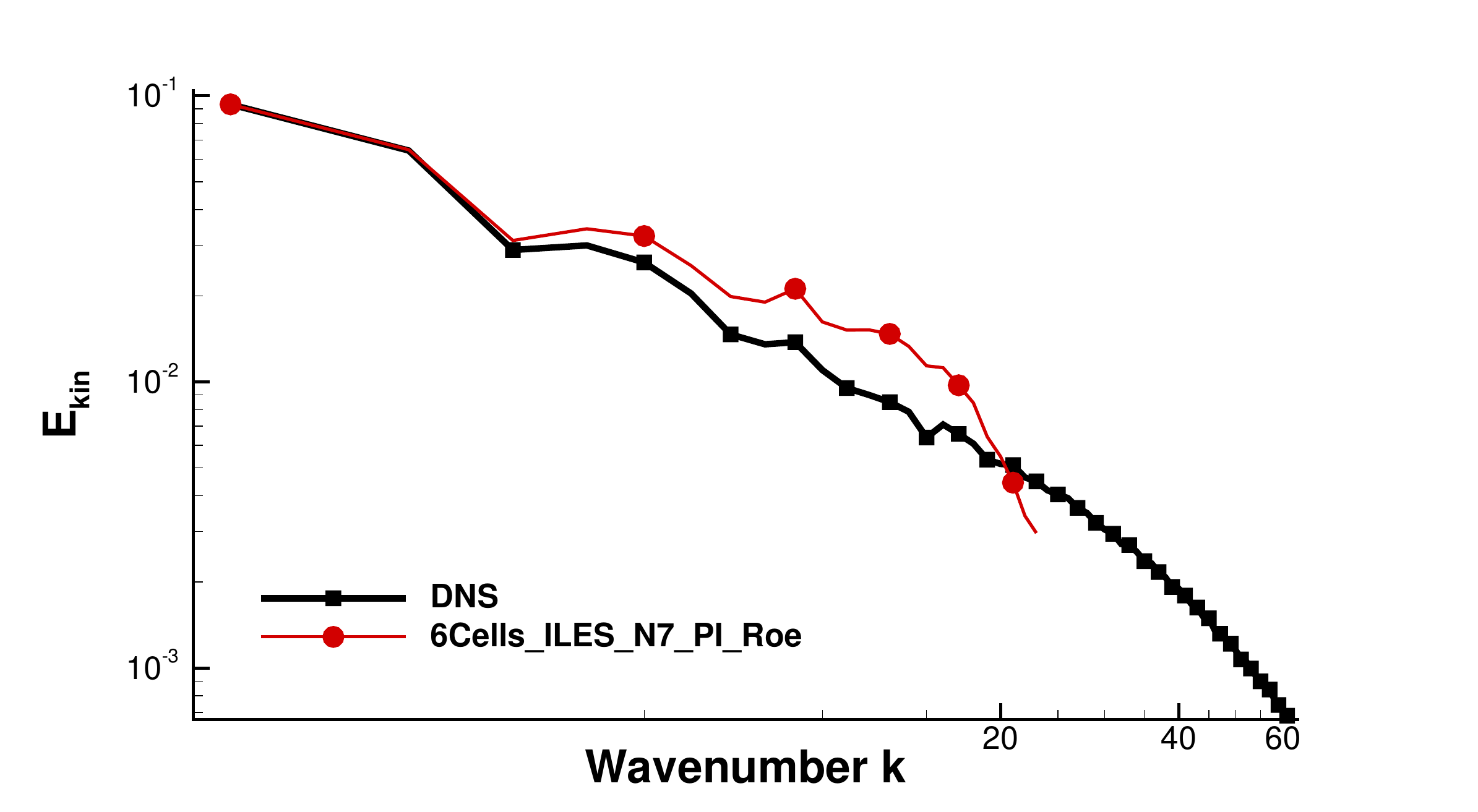}}\\
  \caption{Decaying homogeneous isotropic turbulence DNS/LES kinetic energy spectra: a) Split form DG iLES $144^3$ DOF b) Split form DG iLES $48^3$ DOF. PI interface numerical flux is augmented with additional Roe type matrix dissipation.}
\end{figure}

Figure \ref{fig:iLES12split} shows the result for the case where about $43.2\%$ of dissipation is resolved, i.e. the $144^3$ DOF setup. 

Similar to the standard DG iLES results with polynomial de-aliasing, compare figure \ref{fig:highresOI}, the split form DG iLES produces excellent results for this setup. However, when decreasing 
again the resolution to a more realistic LES, where only $0.9\%$ of the total dissipation is resolved, the iLES DG approach fails again \ref{fig:iLES6split}. 

This is consistent with the DG iLES results obtained with polynomial de-aliasing as discussed in the introduction. We note again that in the iLES approach the dissipation stems solely from the element interface 
contributions. It was shown in \cite{KoprivaGassner_Dispersion,moura2015} that for high order DG schemes, this dissipation is focused on the highest 
wavenumbers. Consequently, this type of dissipation lacks a low wavenumber component, which provides a possible explanation for the failure: the 
theoretical findings of Kraichnan \cite{kraichnan1976eddy} state that for an infinite inertial KE range the resulting eddy viscosity is a constant (plateau) 
in the low wavenumbers and rising to a cusp towards the cut-off wavenumber. The plateau becomes especially relevant for coarse resolutions, or similarly for 
high Reynolds numbers, leading to a cut-off wavenumber in the inertial range. This is discussed in detail by Lesieur e.g. in \cite{lesieur1996}.

\begin{figure}[!htpb]
   \centerline{\includegraphics[width=0.5\textwidth, height=4.5cm, trim=20 10 60 20 ,clip]{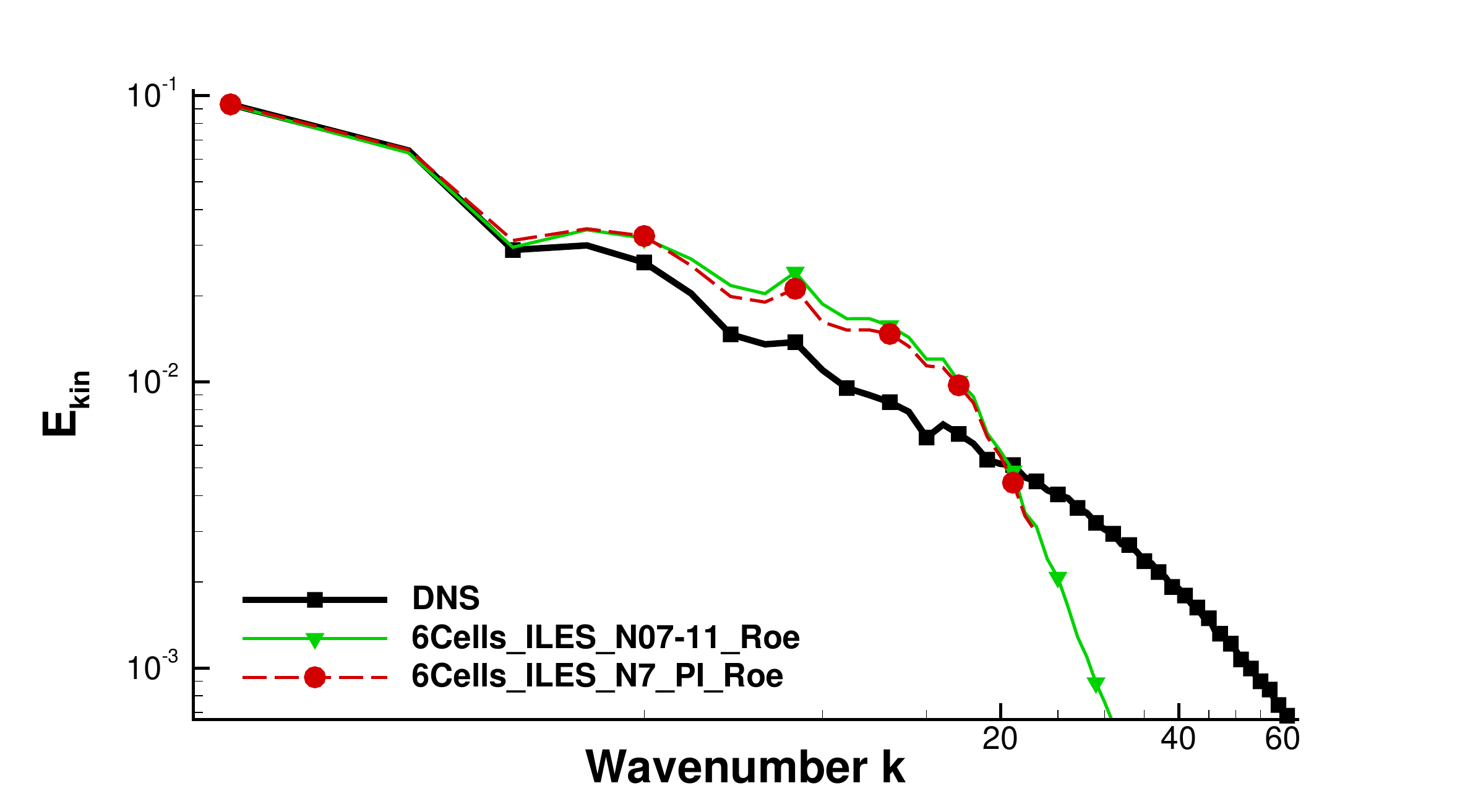}\includegraphics[width=0.5\textwidth, height=4.5cm, trim=20 10 60 20 ,clip]{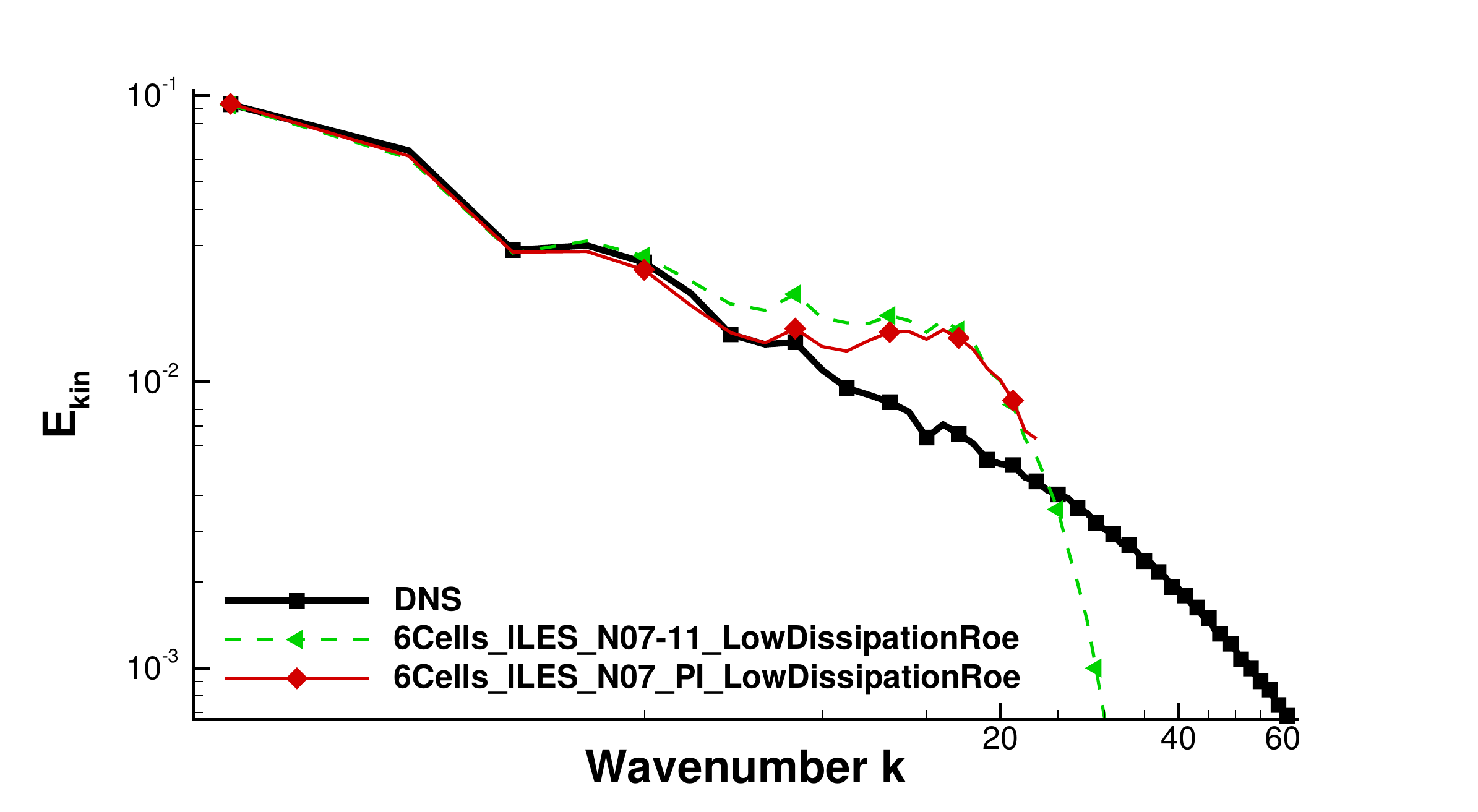}}
   \caption{Decaying homogeneous isotropic turbulence DNS/LES with iLES results based on different de-aliasing strategies: red is based on the PI split form and green uses polynomial de-aliasing.  Both setups use the dissipation term of Roe's Riemann solver at the element interfaces.}
    \label{fig:iLES_dealiasing}
\end{figure}

The left part of Figure \ref{fig:iLES_dealiasing} shows the comparison of the coarse iLES with the two different de-aliasing strategies, namely the polynomial de-aliasing and 
the kinetic energy preserving split form. Both discretizations use the dissipation terms of Roe's Riemann solver at the element interfaces. The almost identical KE spectra show clearly that in this case, the failure 
is due to the unsuited dissipation added by the Riemann solver, and not linked with the de-aliasing procedure. We note that a similar behaviour can be observed with Lax-Friedrichs type dissipation added at element interfaces. 

Finally, we also plot the results obtained with the dissipation terms of the low dissipation Riemann solver of O\ss wald et al. \cite{osswald2015l2roe} in the right part of Figure \ref{fig:iLES_dealiasing}. This variant is a modification to the original Roe solver to account for low Mach number flows (DHIT maximum local Mach number is about $0.1$, on average about $0.02$). This solver gives better results for DHIT when using finite volume schemes \cite{osswald2015l2roe}, compared to the classical Roe solver. The idea of the modification is to multiply the jump in the velocities at the interface by the local Mach number at the interface, resulting in a reduced dissipation in the momentum equations for low Mach number flows. We note that low dissipation Riemann solvers were also successfully used in a high order iLES DG approach with polynomial de-aliasing by \cite{diosady2015}, with an extension to a dynamic variational multi-scale model for the test case of plane channel flow in \cite{murman2016}. The authors clearly demonstrated the need for low dissipation Riemann solvers. 

In contrast to the classic Roe dissipation results in the left part of Figure \ref{fig:iLES_dealiasing}, in the right part of the Figure with the low dissipation Riemann solver of O\ss wald et al. the two de-aliasing strategies behave differently. Indeed one can observe generally a slight improvement of the results for mid range scales, whereas the accuracy for higher wavenumbers is still quite low, underlining again the failure of the iLES approach for very coarse resolutions. It is worth pointing out that split form iLES DG variant is slightly more accurate for mid range scales compared to the polynomial de-aliasing approach and furthermore behaves somewhat similar to the no-dissipation case \ref{fig:split_cent} for the small scales. This is due to the very low dissipation of the low Mach number modification in the momentum equations and thus due to a very low overall kinetic energy dissipation in combination with the kinetic energy preserving split form.  We will demonstrate in section \ref{sec:detailed turbulence modeling} that the split form result with the modified Roe Riemann solver can be drastically improved when an actual SGS turbulence model is added to account for the missing small scale dissipation in the momentum equations instead of using a Riemann solver type dissipation term in the momentum equations.

\section{Split Form DG LES}
\label{sec:kep for expl}

Our second strategy is to use a purely explicit LES DG strategy based on the dissipation free split form result \ref{fig:split_cent}. Instead of adding numerical dissipation at the element interfaces via the Riemann solver, we show in this section the possibility of using explicit turbulence modelling for DHIT to account for the missing small scale dissipation. The results are remarkable, as we will show that by adding the simple Smagorinsky model already gives superior results compared to the iLES DG results discussed above. We note that for dissipation-free base discretisations, the standard Smagorinsky model can be expected to produce reasonably good results if the model constant is adjusted to the underlying discretisation scheme. Up to now, this strategy is only applicable for the novel dissipation free split form DG. 

\begin{figure}[!htpb]
    \centerline{\includegraphics[width=0.5\textwidth, height=4.5cm, trim=20 10 60 20 ,clip]{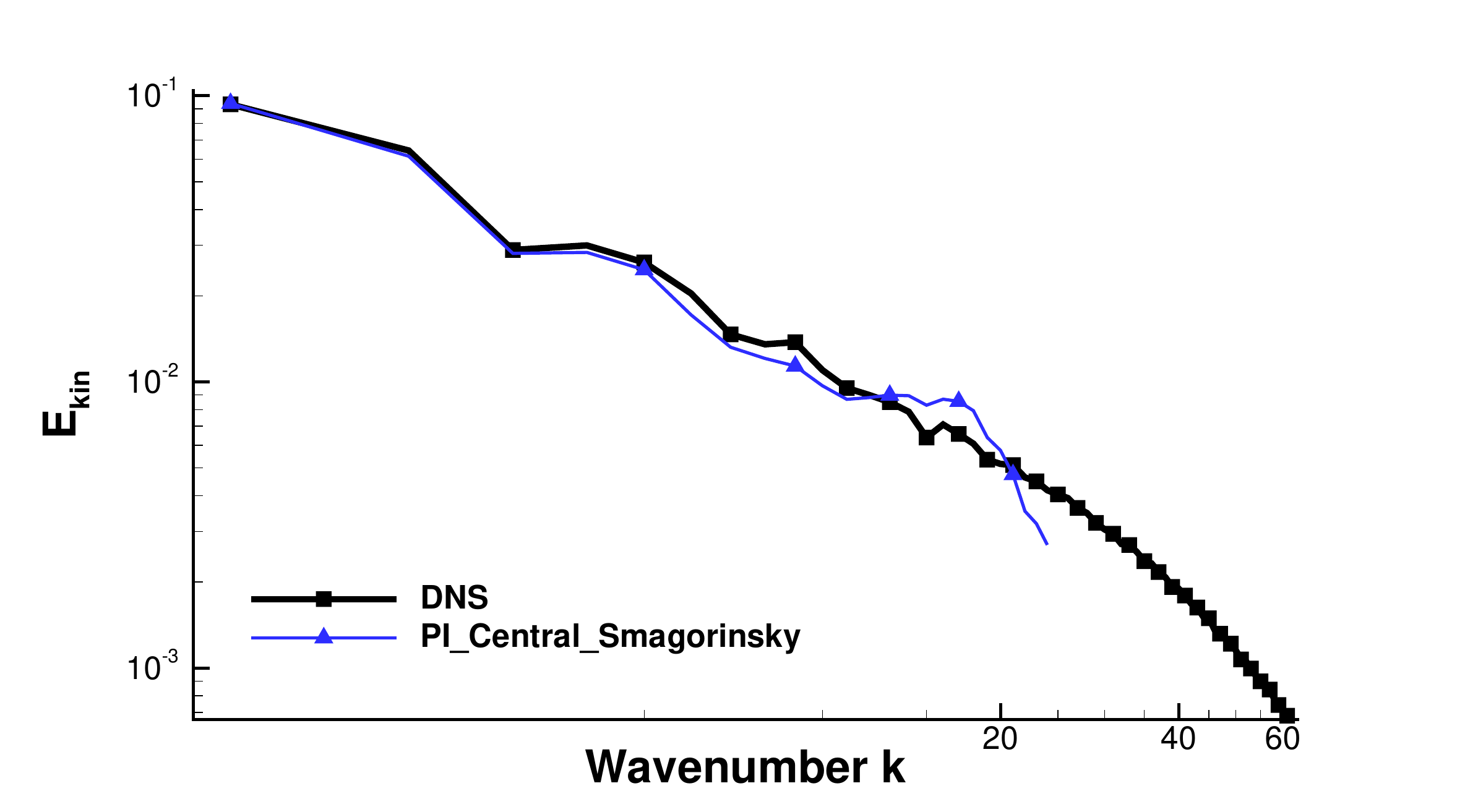}\includegraphics[width=0.5\textwidth, height=4.5cm, trim=20 10 60 20 ,clip]{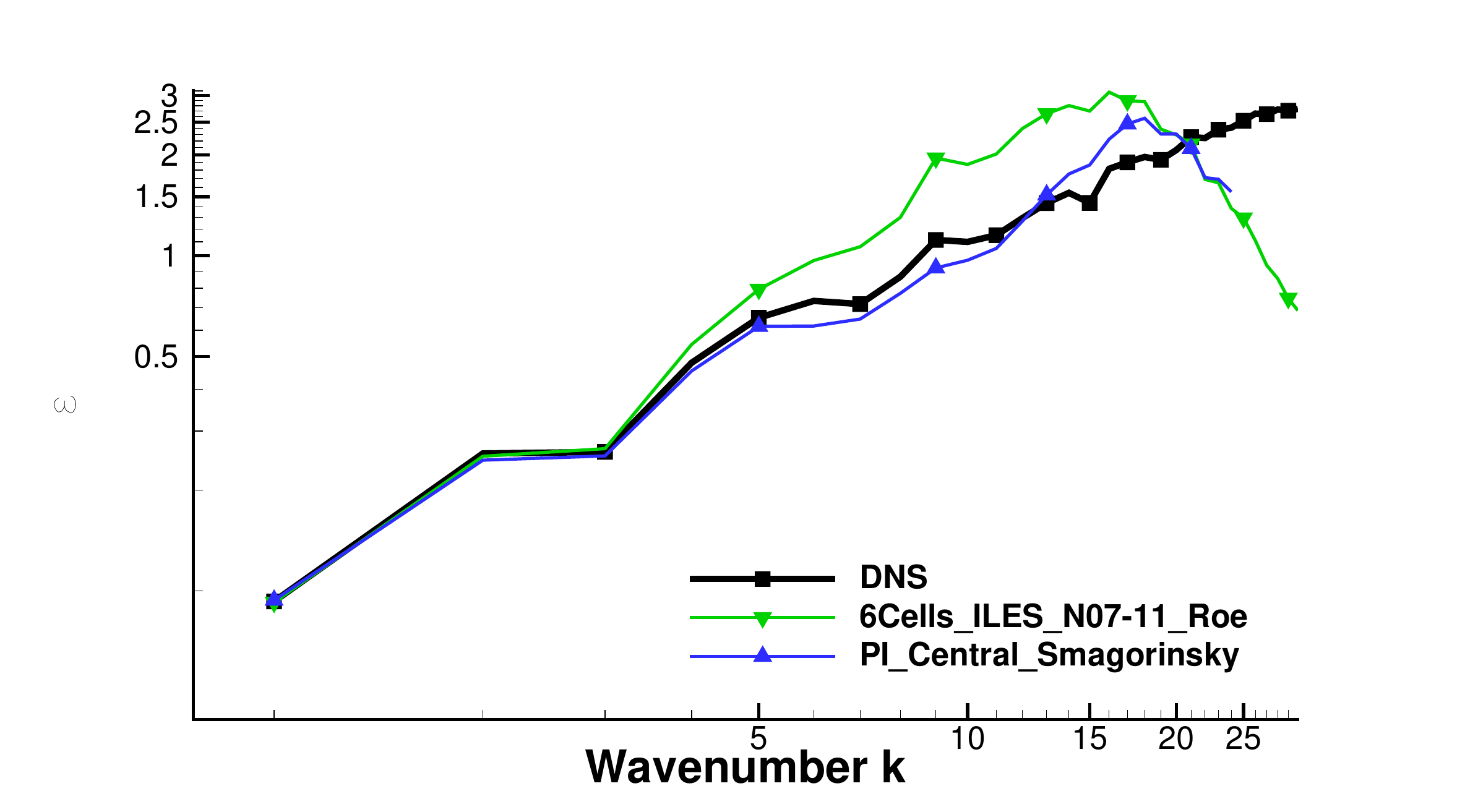}}
      \caption{\label{fig:split_cent_Smago}Left: Decaying HIT KE spectra for the dissipation free baseline scheme with explicit Smagorinsky model with constant $C_s=0.12$. Right: Comparison of dissipation rate spectra to the iLES DG approach with Roe type dissipation and polynomial de-aliasing.}
\end{figure}

In the left part of Figure \ref{fig:split_cent_Smago} the KE spectra is shown, adding only Smagorinsky's model \cite{smago} to the dissipation free baseline scheme. The result is already acceptable, and for the first time a successful coarse resolution LES DG result is obtained. 

We further plot in the right part of Figure \ref{fig:split_cent_Smago} a comparison of the spectra of the dissipation rate $\epsilon$ to the iLES DG approach with classic Roe type dissipation and polynomial de-aliasing. It is remarkable that even with a rather simple and limited SGS model, the result of the explicit LES DG is clearly superior to the iLES DG result. The key to this result is a baseline scheme that is virtually dissipation free and at the same time robust for under-resolved simulations. Thus, the SGS model dissipation is not necessary to retain stability by adding dissipation to the dissipation free baseline scheme, but is used to account solely for missing small scale dissipation. This strict separation of stability and model dissipation also allows for a more precise analysis of the affect of the turbulence model: The spectra shows a slight energy underestimation of kinetic energy in the medium range wavenumbers. This is usually a secondary effect due to undamped small scale fluctuations. This modelling effect is well known and various examples can be found in literature, an overview is given again by e.g. Lesieur \cite{lesieur1996}. 

Nevertheless, the dissipation-free baseline scheme with explicit Smagorinsky model offers a superior result compared to the state of the art iLES DG strategy. It is thus reasonable to proceed and investigate more advanced turbulence models and assess if in contrast to the iLES approach advanced turbulence modelling within the split form DG approach actually improves the simulation. 

\section{Advanced Turbulence Modelling for Split Form DG LES}
\label{sec:detailed turbulence modeling}

Our third strategy is to add a cusp like eddy viscosity mechanism to cure the secondary effect of undamped small scale fluctuations observed for the standard Smagorinsky approach observed in the previous section: In order to improve the method described in section \ref{sec:kep for expl} we show three possible ways of introducing additional damping for the small scales (cusp). The aim is to show that the dissipation free baseline scheme offers a lot of freedom and possibilities to design and control the overall amount and snape of the dissipation. All strategies are focused on improving the modelling of the missing small scale dissipation, as the aspect of robustness and numerical stability is already accounted for by using the dissipation free baseline split form DG scheme. The three strategies considered and tested are: 

\begin{itemize}
 \item[i)] adding a low dissipation Riemann solver \cite{osswald2015l2roe} at the interface,
 \item[ii)] adding a penalty in the viscous numerical flux (BR2 with $\eta_{BR}=2$),
 \item[iii)] adding additional eddy-viscosity in a variational multi scale sense with a cusp profile over polynomial modes 
following Chollet et al. \cite{chollet1985}.
\end{itemize}

The first two methods implicitly add dissipation via the numerical interface fluxes in addition to the Smagorinsky model. In principal, these can be seen as a mixture of iLES and explicit LES. However, as pointed out above, the interface dissipation acts like a high frequency filter and affects the small scales while the influence on the rest of the scales is very small, i.e. it introduces an additional damping of the small scales only (cusp). We note that in DG it is possible to introduce numerical interface dissipation via the advective Riemann solver or the viscous numerical flux function. For the advective numerical dissipation, we chose the low dissipation Roe variant discussed above, as the classic Roe dissipation already destroys the accuracy as shown in previous sections. For the numerical dissipation of the viscous terms, we extend change the BR1 approach to the BR2 approach, where the an additional penalty term of the jump of the solution is added to the local gradients. The effect of this penalty can be adjusted by the BR2 constant $\eta_{BR}$. In contrast to the first two strategies, method iii) can be seen as a purely explicit LES with an advanced SGS model where dissipation is increased for the highest polynomial modes, which correspond to the small scale contents of the solution. 

\begin{figure}[!htpb]
   \centerline{\includegraphics[width=0.7\textwidth, height=6.3cm, trim=20 10 3 8 ,clip]{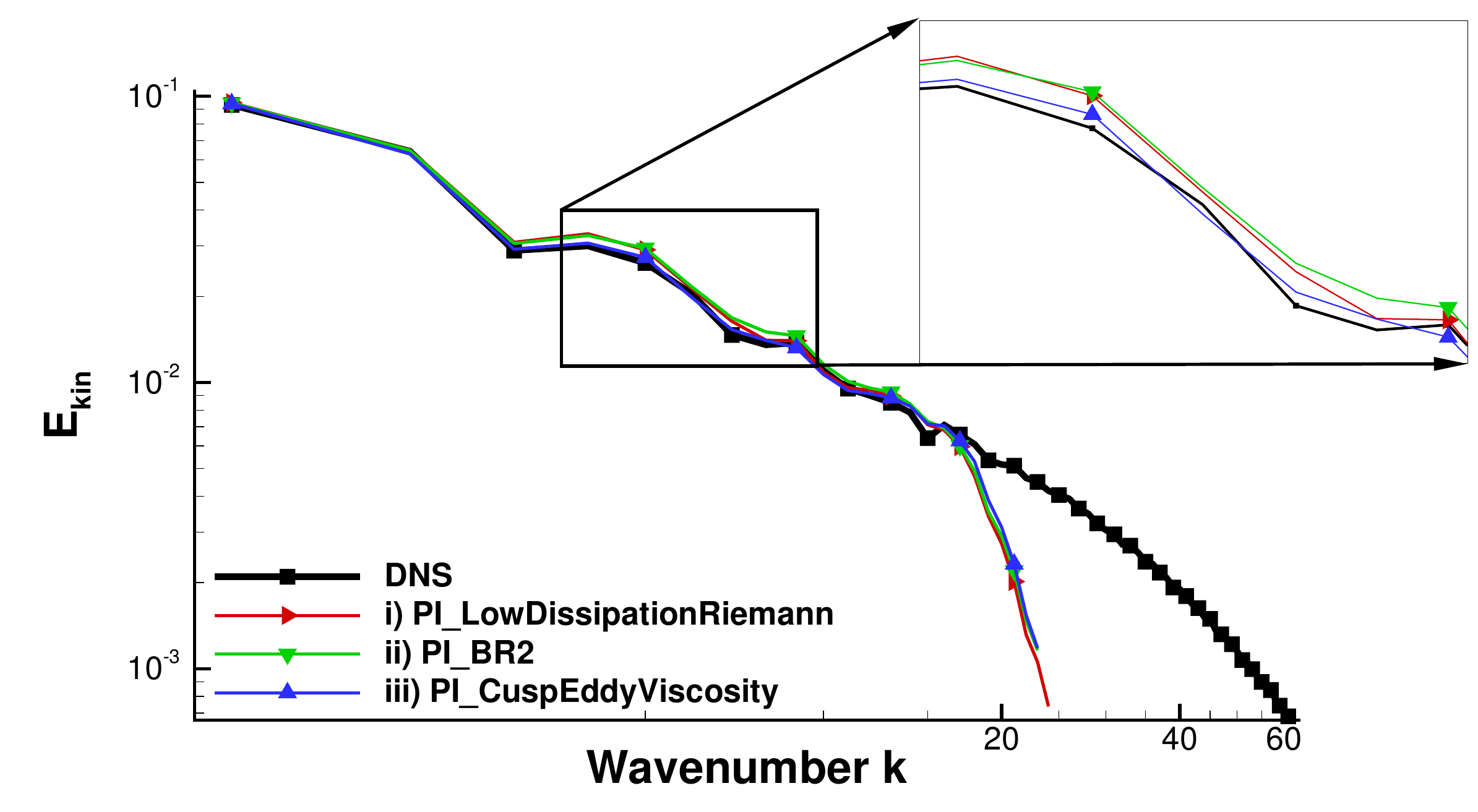}}
   \caption{Decaying homogeneous isotropic turbulence DNS/LES with different modelling approaches. Results for modelling approaches almost collapse. Zoom in to 
region with largest differences.}
    \label{fig:split_lesGood}
\end{figure}

The results of the three advanced strategies are shown in Figure \ref{fig:split_lesGood} and all three variants give an improved KE spectrum. The results of our investigations are that method iii) allows for the best fit with the DNS, if the model parameters are adjusted to the dissipation free split form DG baseline scheme. Details of this model are given in \ref{App:appendixPlCusp}.

\begin{figure}[!htpb]
   \centerline{\includegraphics[width=0.5\textwidth, height=4.5cm, trim=2 10 3 8 ,clip]{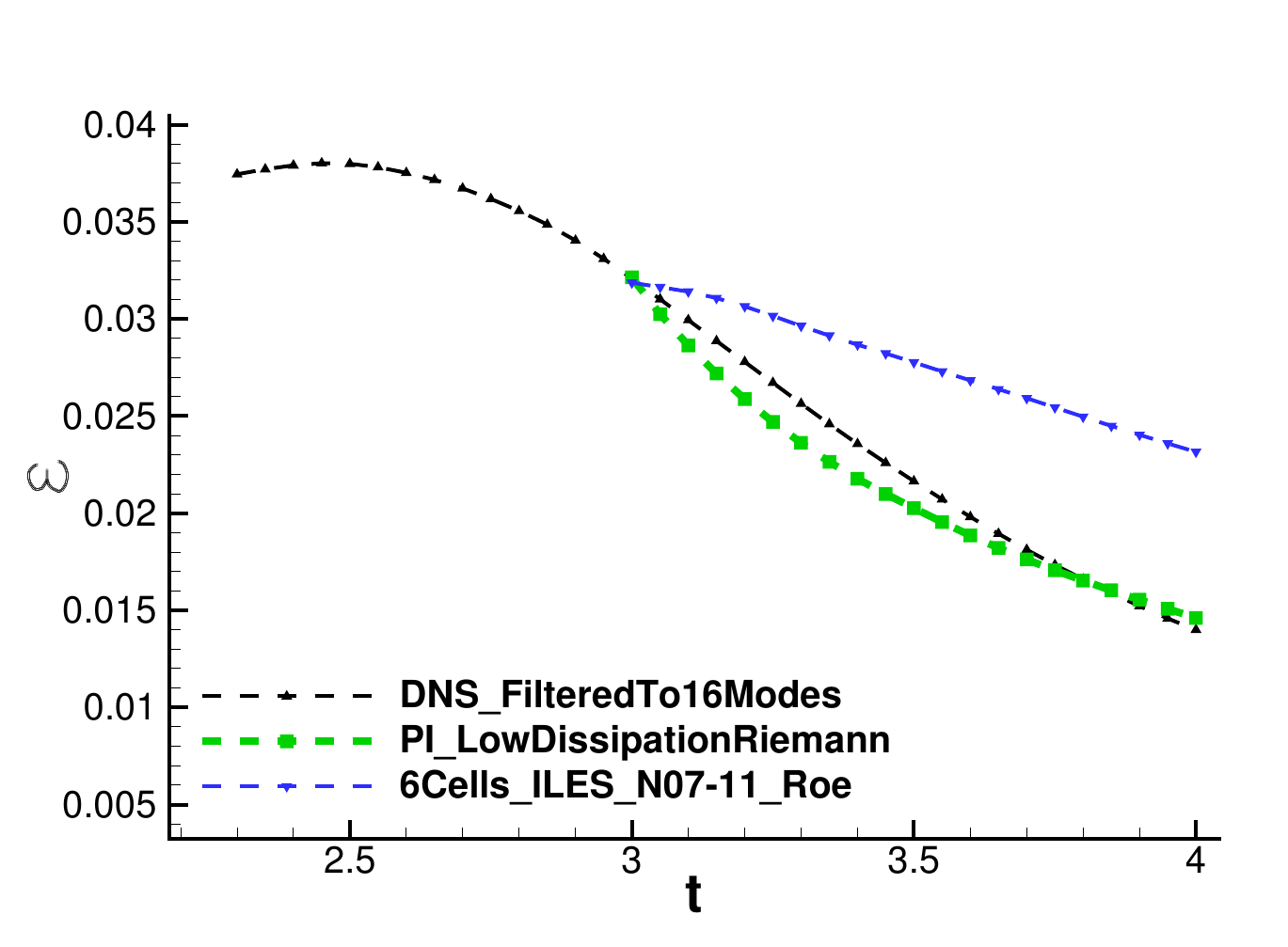}\includegraphics[width=0.5\textwidth, height=4.5cm, trim=2 10 3 8 ,clip]{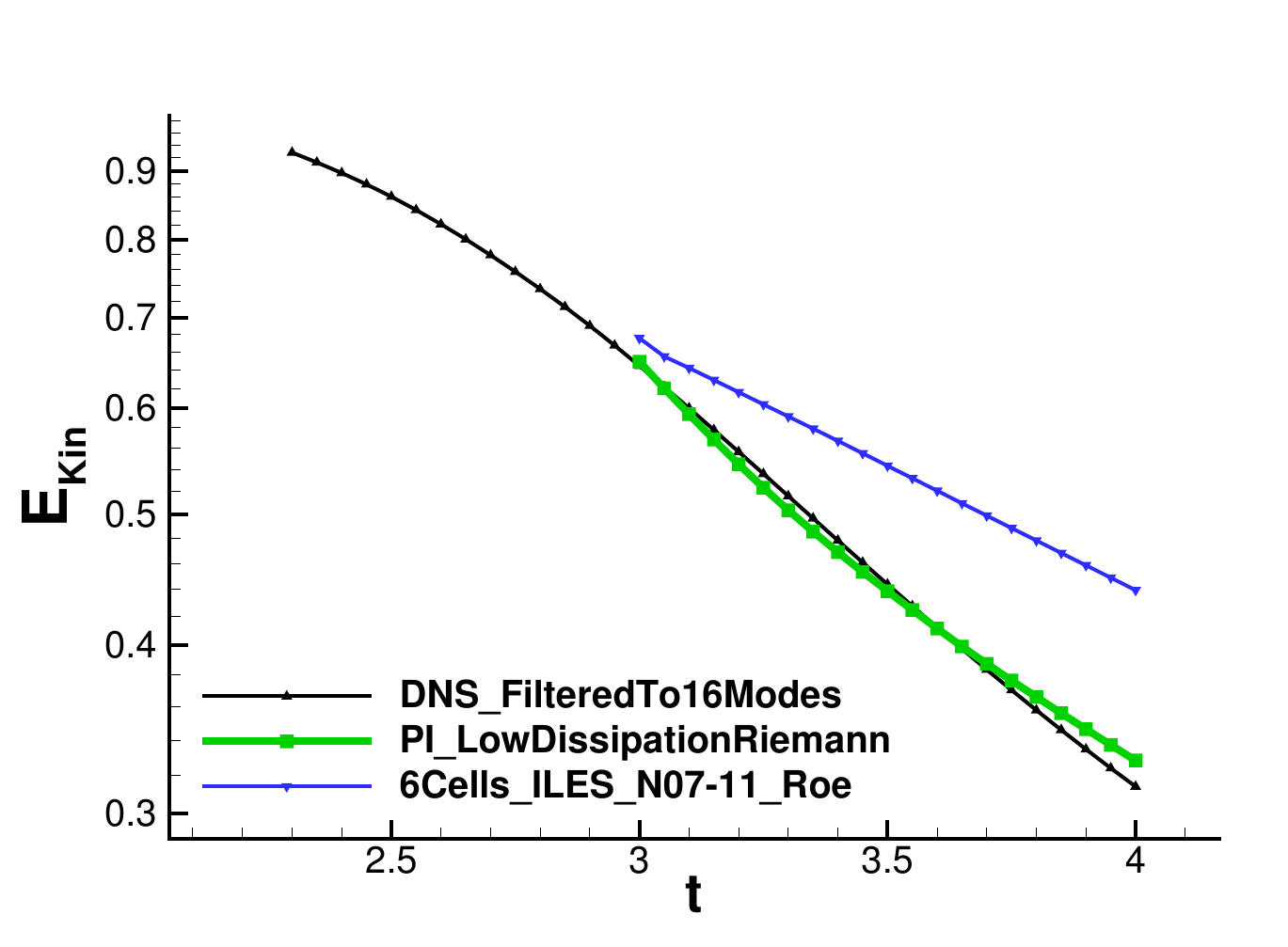}}
   \caption{Comparison of dissipation rate (left) and kinetic energy (right) over time for (filtered DNS, split form LES DG, and standard iLES DG. Dissipation rate and kinetic energy are integrated in 
Fourier space up to $k=16$.}
    \label{fig:diss_overTime}
\end{figure}

Focusing on strategy i) with the low dissipation Roe Riemann solver, Figure \ref{fig:diss_overTime} shows a comparison to the standard iLES DG approach with Roe dissipation and polynomial de-aliasing of the dissipation rate and the kinetic energy over time. The dissipation (kinetic energy) rate is computed in Fourier space as $\epsilon = 2\nu\int_{k=1}^{16}k^2E(k)dk$ ($E_{Kin} = 
\int_{k=1}^{16}E(k)dk$) to obtain comparable results between LES and (filtered) DNS. We find that the standard iLES DG approach is not capable of 
reproducing the energy decay over time with this coarse resolution. In comparison, the novel method i) is in much better agreement with the (filtered) DNS result. 

\subsection{Pressure-dilatation fluctuations}

\begin{figure}[!htpb]
  \subfigure[Central+Smagorinsky (section \ref{sec:kep for expl})]{%
    \label{fig:diss cent}
    \includegraphics[width=0.5\linewidth, trim=10 20 140 150,clip]{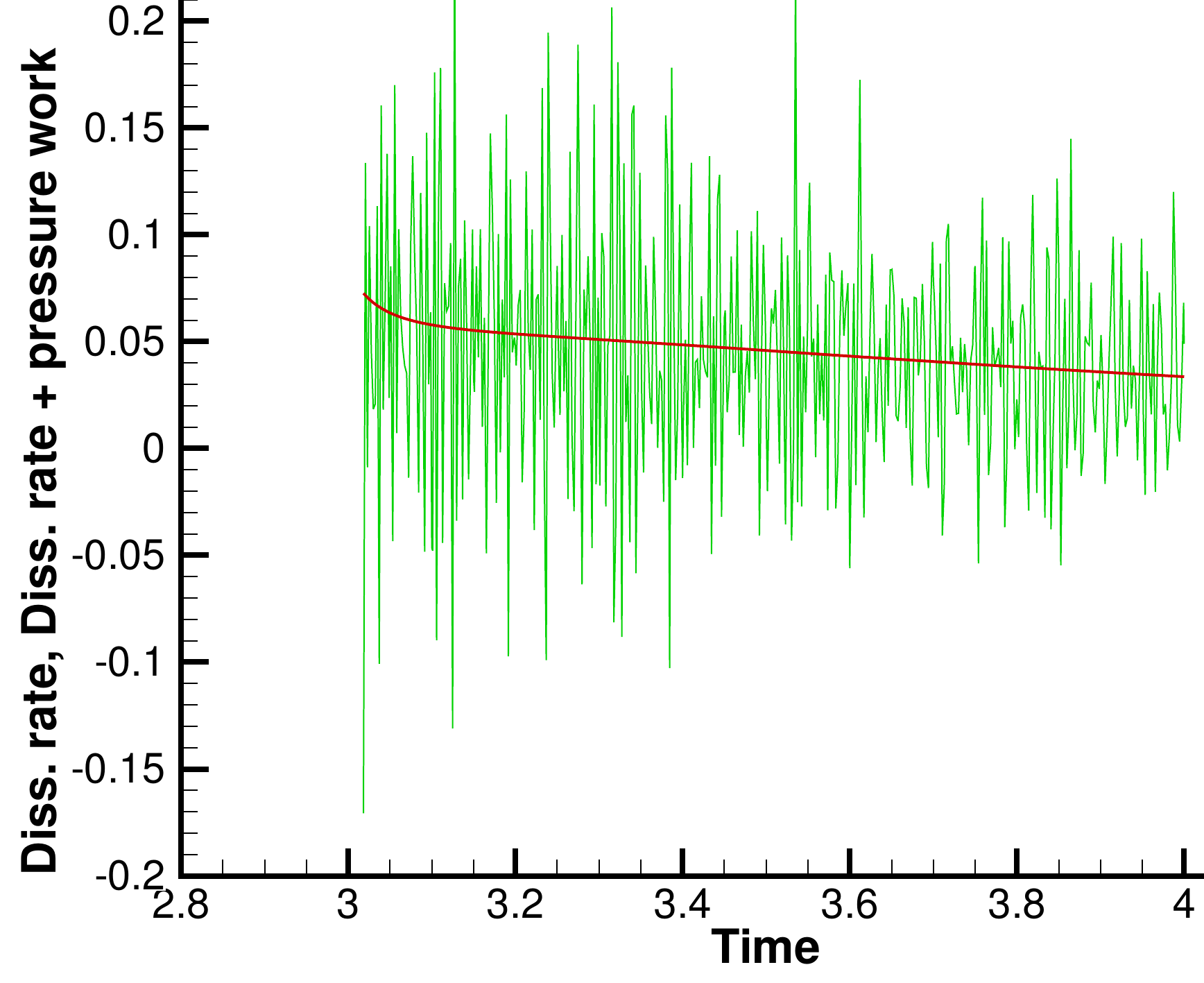}}
  \subfigure[Low dissipation Roe, method i) (section \ref{sec:detailed turbulence modeling})]{%
    \label{fig:diss roel2}
    \includegraphics[width=0.5\linewidth, trim=10 20 140 150,clip]{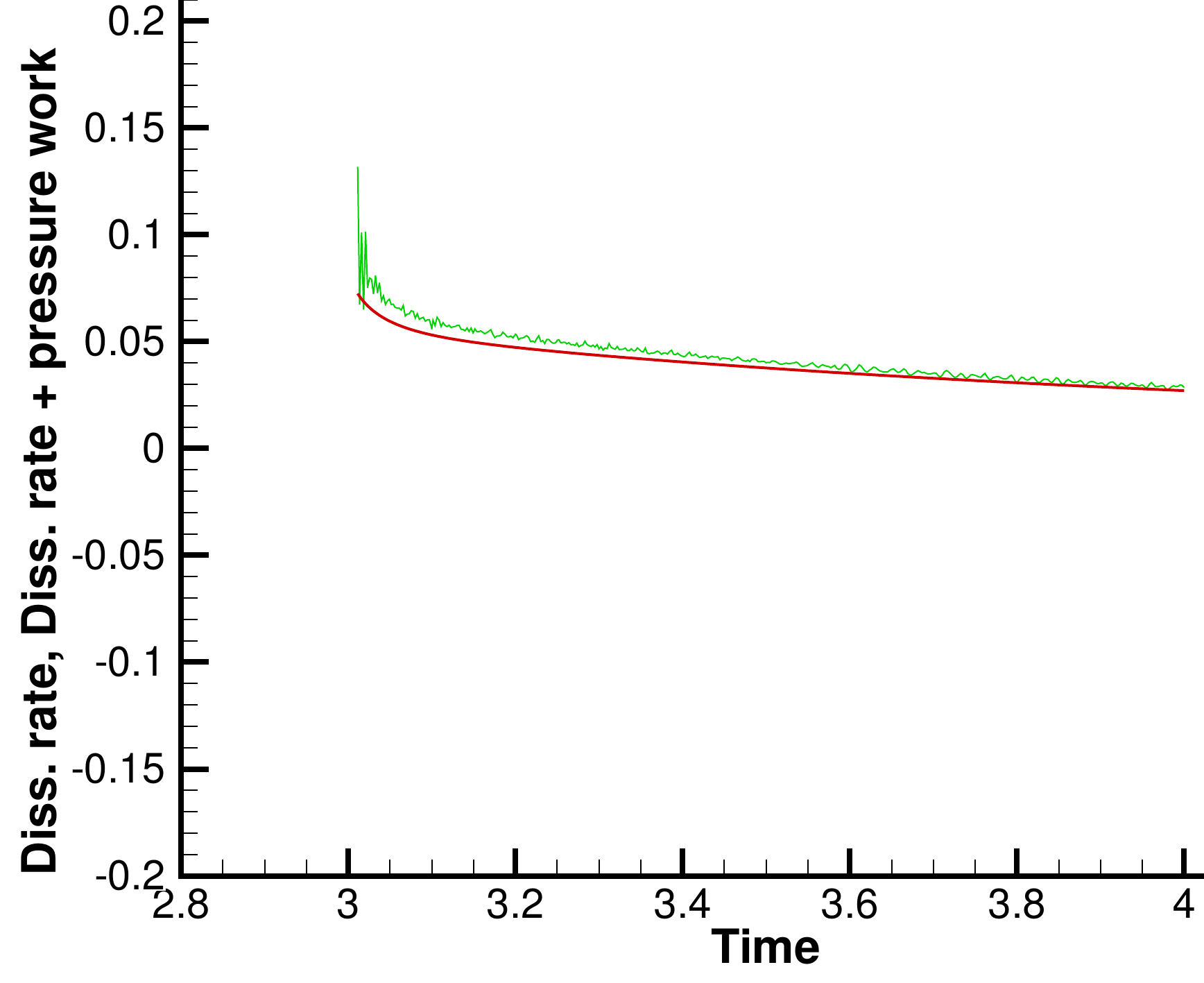}}\\
  \subfigure[BR2, method ii) (section \ref{sec:detailed turbulence modeling})]{%
    \label{fig:diss br2}
    \includegraphics[width=0.5\linewidth, trim=10 20 140 150,clip]{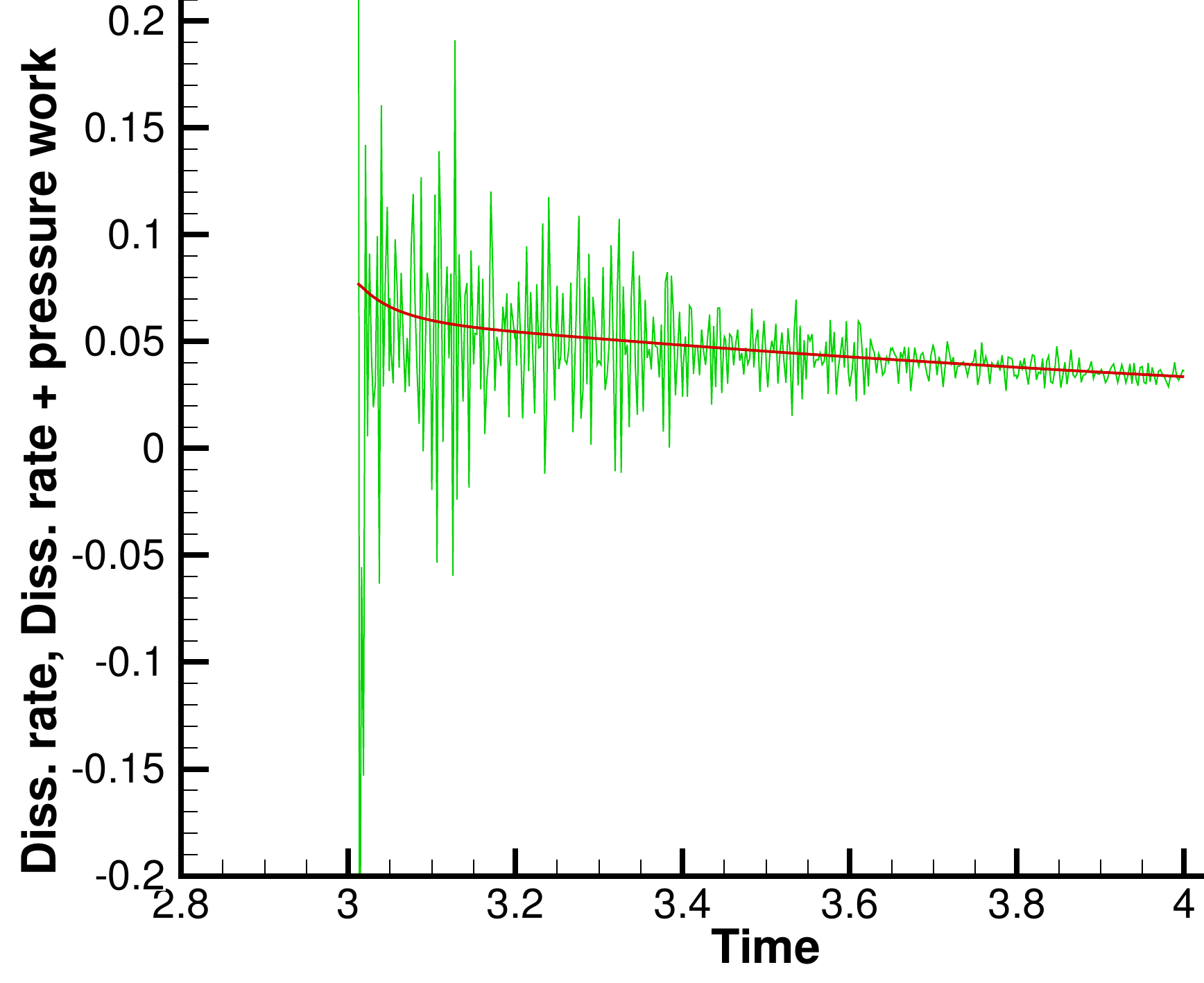}}
  \subfigure[Plateau-Cusp, method iii) (section \ref{sec:detailed turbulence modeling})]{%
    \label{fig:diss plcusp}
    \includegraphics[width=0.5\linewidth, trim=10 20 140 150,clip]{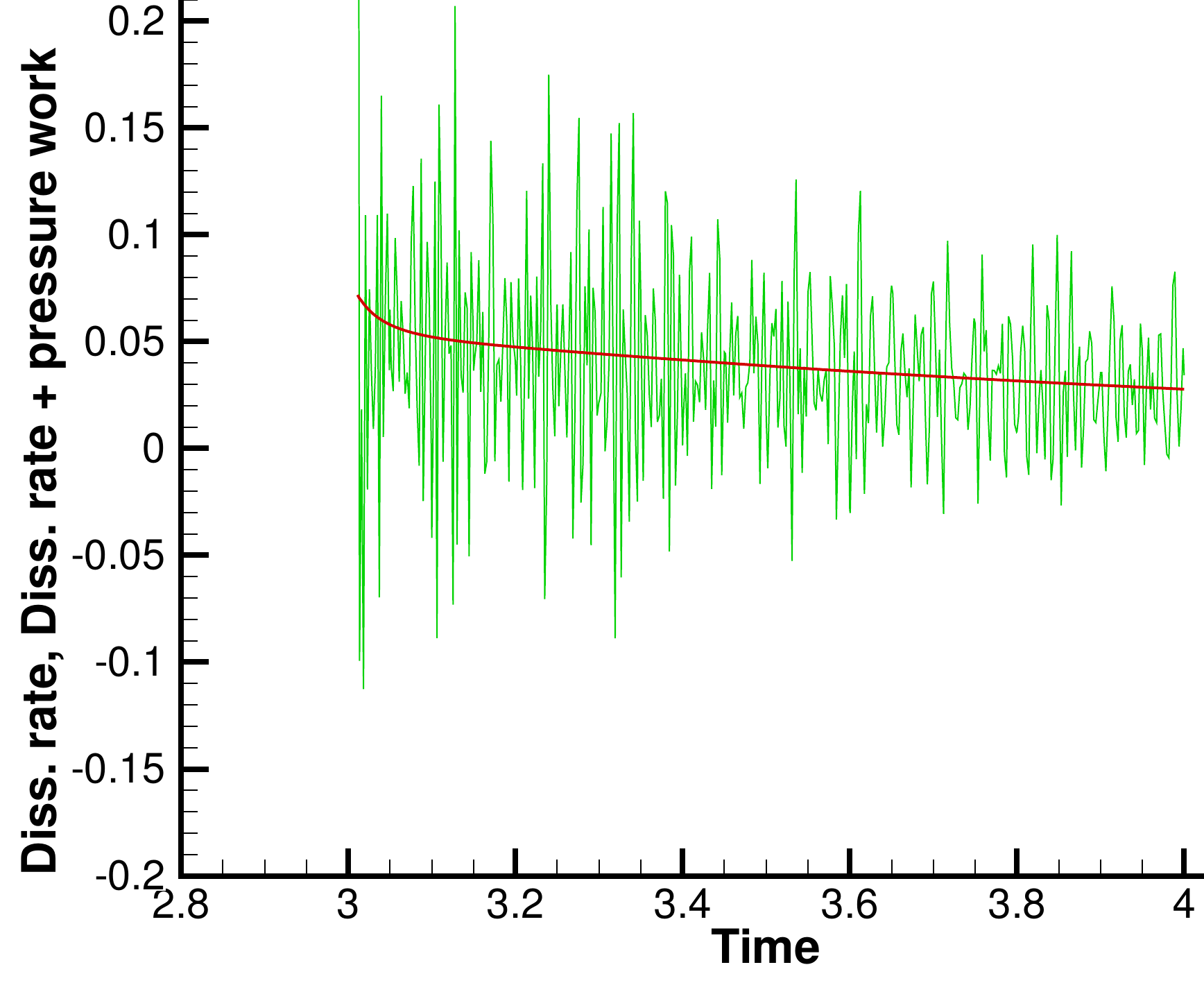}}
 \caption{ \label{fig:diss}Decaying homogeneous isotropic turbulence LES dissipation rate over time for different modelling approaches. Red: 
$\epsilon=2\nu\int_{\Omega}S_{ij}S_{ij}d\Omega$, green: $\epsilon+\int_{\Omega}(\nabla \circ \mathbf{U})pd\Omega$}
\end{figure}

An interesting side-effect of the advanced modelling results presented above is that it is now possible to observe other, formerly hidden, artefacts in the LES solution: pressure work can oscillate and strong fluctuations can be observed. It is important to note that by construction, neither the kinetic energy preserving PI split form nor the standard Smagorinsky model of the standard explicit LES DG approach presented above provide mechanisms to damp pressure of density fluctuations, especially regarding the inter element discontinuities of the DG approach. Thus, these quantities can and do oscillate as can be seen in the pressure work (syn. pressure-dilatation) $(\nabla \circ \mathbf{U})p$, Figure \ref{fig:diss cent}. 

As investigated by \cite{chandrashekar2012kinetic}, tailored Riemann solver dissipation, damping only density and pressure fluctuations, is sufficient to remove the oscillations. Among the three methods of section \ref{sec:detailed turbulence modeling} only method i), the low dissipation Riemann solver, provides a mechanism to damp these fluctuations, compare Figure \ref{fig:diss}. The fluctuations in method ii) are somewhat smaller due to the added penalty terms in the computations of the gradients of the density and energy, Figure \ref{fig:diss br2}. Method iii) does not provide additional damping of pressure or density fluctuations by construction and thus the pressure work fluctuations in Figure \ref{fig:diss plcusp} remain almost unchanged compared to \ref{fig:diss cent}.

In summary, the best LES results are obtained with method i): the baseline dissipation free split form DG scheme with a low dissipation Riemann solver that controls the fluctuations in density and pressure (but unaffected the dissipation in the momentum equations) with a standard Smagorinsky SGS model that provides suitable dissipation of the momentum equations and thus the kinetic energy.  This clearly demonstrates the freedom in designing and precisely controlling the amount and shape of the added dissipation when the novel split form DG approach is used as a baseline scheme.



\section{Discussion}
\label{sec:conclusion}

\subsection{Summary}

This work presents the first extension of the novel split form DG framework \cite{gassner2016split} to explicit LES and iLES. We have first shown that there seems 
to be a common upper limit in the applicability of DG schemes for implicit LES. Only in cases where a substantial part of the dissipation is resolved by the DG 
approximation, the results are in good agreement with the reference. However, for realistic coarse resolutions, we show that the Riemann solver dissipation in 
combination with either de-aliasing strategy, polynomial de-aliasing or split forms, is not well suited as a sub-grid scale model for turbulence. 

As a remedy, we propose a novel approach based on split forms that give non-linear robustness in combination with kinetic energy preservation, namely the split form introduced by Pirozzoli as a variant of the Kennedy and Gruber split form. This approach allows us to construct a baseline DG scheme with virtually no kinetic energy dissipation that is still stable for vortical dominated turbulence, even in case of severe 
under resolution. The main idea is that this stable but dissipation free DG scheme allows for precise control of the amount and shape of the dissipation added to model the effect of missing small scales, as no additional dissipation is needed to restore the robustness of the scheme in case of under resolution. 

As a first result, this dissipation free baseline scheme in combination with a simple Smagorinsky model gives superior results compared to iLES DG. We further advance these results and introduce the following three strategies to add a cusp like behaviour to increase dissipation of the smallest scales:
\begin{itemize}
 \item a kinetic energy consistent split-DG scheme,
 \item the low dissipation Roe-type Riemann solver of \cite{osswald2015l2roe},
 \item and Smagorinsky's SGS model.
\end{itemize}
All three advanced strategies yield improved results compared to the standard Smagorinsky model. But only method i) effectively controls oscillations in the pressure-dilatation and is thus the clear winner overall. It is also worth noting that method i) substantially improves the results in comparison to start of the art implicit LES DG. The results clearly underline the main idea that it is indeed possible to fully control and precisely adjust the dissipation in split form DG to account for missing small scales. As a final remark, the novel split form LES DG approach considerable reduces the computational cost by at least a factor of $3$ in our tests within our in-house DG simulation framework.

\subsection{Outlook}

\begin{figure}[!htpb]
    \centering
    \includegraphics[width=0.8\textwidth, trim=20 10 50 20 ,clip]{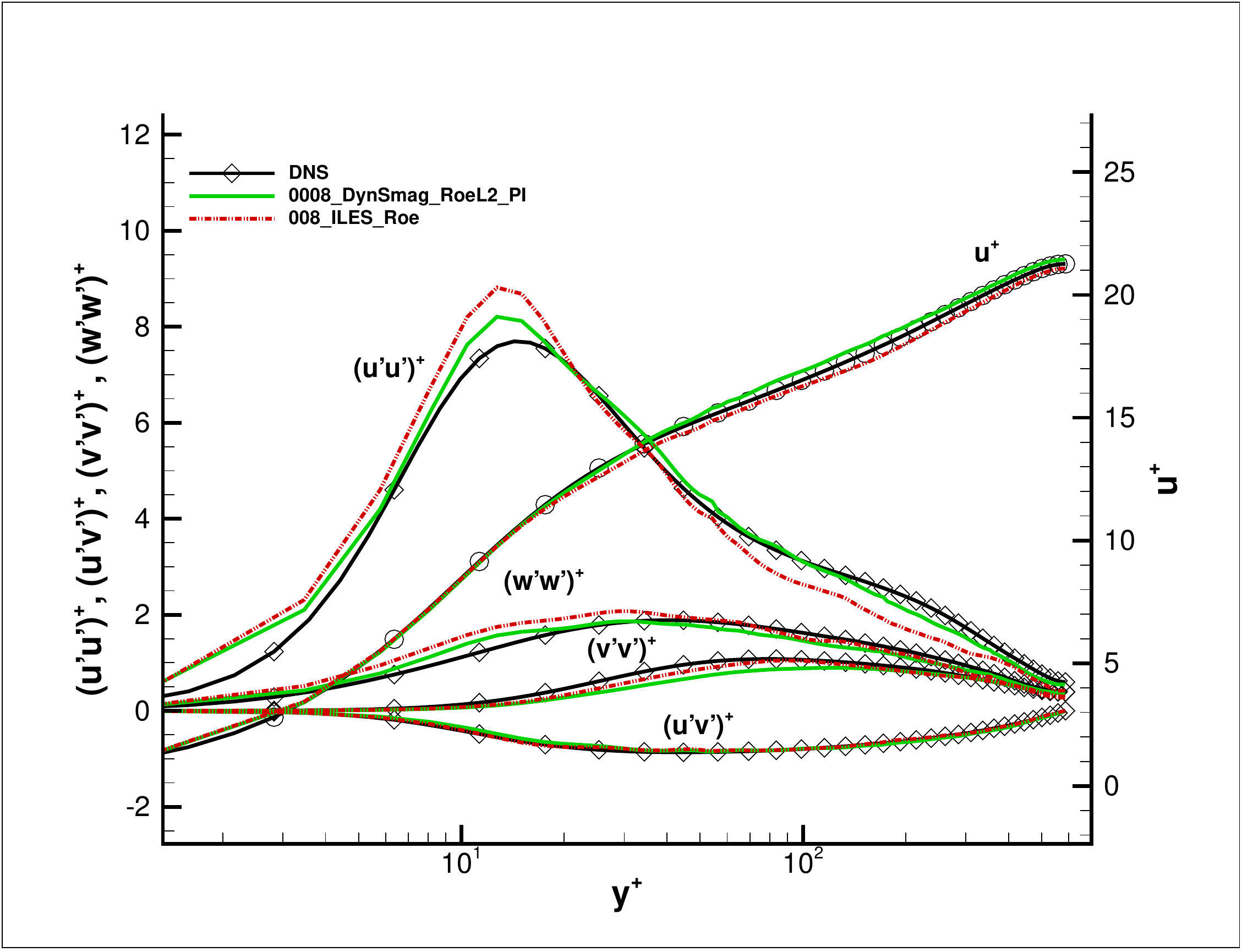}
  \caption{Plane turbulent channel flow, $Re_{\tau}=590$, $8$ cells per direction, $N=7$. Green line shows the novel split form DG LES approach with an explicit dynamic Smagorinsky model. Red curve shows standard implicit LES approach with polynomial de-aliasing and Roe's Riemann solver.  Reference DNS is taken from Moser \cite{moser1999}, 
and no model LES}
  \label{fig:channel}
\end{figure}

The split form DG LES approach can be directly applied with more complex SGS models for more demanding applications. As a first result we show the computation of a plane turbulent channel flow employing the dynamic Smagorinsky model of Germano \cite{germano91}, with Lilly's modification \cite{lilly92}. The flow is simulated for $Re_{\tau}=590$, with a $8^3$ grid and $N=7$ in a domain $[x,y,z]\in[0,2\pi]\times[-1,1]\times[0,\pi]$. Grid spacings are $\Delta x^+\approx58$, $\Delta z^+\approx29$ and $\Delta y^+_{min/max}\approx7/28$, based on an equidistant inner cell point distribution. The result is compared to the DNS of \cite{moser1999} and the standard iLES DG approach with polynomial de-aliasing and Roe's Riemann solver in Figure \ref{fig:channel}.

Both LES strategies  yield acceptable results, even with this quite coarse resolution. However, the novel split form DG LES approach with explicit SGS turbulence modelling gives better Reynolds 
stress profiles, showing again excellent agreement with the DNS.

\section{Acknowledgement}

G.G has been supported by the European Research Council (ERC) under the European Union's Eights Framework Program Horizon 2020 with the research project \textit{Extreme}, ERC grant agreement no. 714487.

\newpage
\appendix
\section{Smagorinsky model}
\label{App:smago}
The classical Smagorinsky SGS stress, is defined as
\begin{equation}
 \tau_{ij} = -2C_S^2\Delta^2|S(\nabla U)|(S_{ij}(\nabla U)-\frac{\delta_{ij}}{3}S_{ii}(\nabla U))
\end{equation}
Where for the velocity gradients we use the lifted gradients. The filter width $\Delta$ is computed as the local cell volume divided be $(N+1)^3$. As usual $S$ 
is the strain rate tensor.
\section{High pass filtered Smagorinsky (VMS)}
\label{App:VMS}
The high pass filtered Smagorinsky model (or variational multiscale Smagorinsky model, differing in the fiter operation only. We use a Galerkin projection as 
filter, i.e. a variational multi scale formulation) became popular recenlty, especially for wall bounded turbulence as the viscosity does not affect the 
large scales. The eddy viscosity is computed from the small scale gradients. 
We use a polynomial projection high-pass filter and apply it to the lifted gradients, denoting those as $\widetilde{\nabla U}$, cutting modes grater than 
$N_{filter}$. 
We compute the SGS stress as:
\begin{equation}
 \tau_{ij}^{VMS} = -2C_S^2\Delta^2|S(\widetilde{\nabla U})|(S_{ij}(\widetilde{\nabla U})-\frac{\delta_{ij}}{3}S_{ii}(\widetilde{\nabla U}))
\end{equation}
\section{Plataeu-Cusp model}
\label{App:appendixPlCusp}
In spectral methods it has proven advantageous to use an eddy viscosity with a plateau in the large scales and an a cusp towards the cut of wavenumber 
\cite{chollet1985}. This is usually called the spectral eddy viscosity model (SPEVM). The shape is given in Fourier space by:
\begin{align}
 \nu^t(k|k_c)=\nu^{t+}(k|k_c)\left(\frac{E(k_c,t)}{k_c}\right)^{(1/2)}\\
 \nu^{t+}(k|k_c)=0.267(1+34.5e^{-3.03(k_c/k)})
\end{align}
We seek to add a similar cusp in modal polynomial space. We reformulate the viscosity distribution as:
\begin{equation}
 \nu^{t+}(N+1|N_c+1)=9.21e^{-3.03((N_c+1)/(N+1))})
\end{equation}
This distribution we use as input to a polynomial projection high-pass filter and apply it to the lifted gradients, denoting those as $\widetilde{\nabla U}$. 
We compute the SGS stress in analogy to the high pass filtered eddy viscosity models (or VMS) as:
\begin{equation}
 \tau_{ij}^{cusp} = -2C_{cusp}\Delta^2|S(\widetilde{\nabla U})|(S_{ij}(\widetilde{\nabla U})-\frac{\delta_{ij}}{3}S_{ii}(\widetilde{\nabla U}))
\end{equation}
This ''cusp'' stress we add to the classical Smagorinsky SGS stress, defined equivalent without high-pass filtered gradients
\begin{equation}
 \tau_{ij} = -2C_S\Delta^2|S(\nabla U)|(S_{ij}(\nabla U)-\frac{\delta_{ij}}{3}S_{ii}(\nabla U))
\end{equation}
For the result above we chose $C_{Cusp}=0.33^2$ and $C_S=0.13^2$.

\section{Dynamic Smagorinsky}
\label{App:dynSmag}
The SGS stress of the classical Smagorinsky model is defined as
\begin{equation}
 \tau_{ij} = -2C_S\Delta^2|S(\nabla U)|(S_{ij}(\nabla U)-\frac{\delta_{ij}}{3}S_{ii}(\nabla U)).
\end{equation}
The dynamic procedure determines the constant $C_S$ depending on the flow:
\begin{equation}
 C_S \Delta^2 = -\frac{1}{2}\frac{<L_{ij}M_{ij}>}{<M_{ij}M_{ij}>}
\end{equation}
with $<>$ an average operator over homogeneous directions. We chose to average over Y-Planes within a DG cell.
The Germano identity is defined as
\begin{equation}
 L{ij} = \widetilde{u_iu_j} - \widetilde{u_i}\widetilde{u_j}
\end{equation}
with $\tilde{.}$ denoting a Galerkin projection filter operation, applied in x-z planes. We chose the test filter to $N_{testfilter}=3$, while the 
computational polynomial degree was 
$N=7$. That results to a filter width ratio of $(N)/(N_{testfilter})\approx2.3$.\\
Also
\begin{equation}
 M_{ij} = \widetilde{|S|(S_{ij}(\nabla U)-\frac{\delta_{ij}}{3}S_{ii}(\nabla U))} - \frac{\widetilde{\Delta}}{\Delta}|S(\widetilde{\nabla 
U})|(S_{ij}(\widetilde{\nabla U})-\frac{\delta_{ij}}{3}S_{ii}(\widetilde{\nabla U}))
\end{equation}

\section{Channel flow details}
\label{App:channel}
For brevity we omitted details of the channel flow setup above. In this appendix we complement this. For the channel flow the second method of Bassi and Rebay 
(BR2) \cite{bassi2} was used. The penalty constant $\eta_{BR2}$ is chosen to $1$ at all inner (i.e. not wall boundaries) cell interfaces. In contrast to the 
first method of Bassi and Rebay (BR1), used for the DHIT simulations, the BR2 scheme is not artificial dissipation free. By choosing $\eta_{BR2}=1$ the 
methods are similar in that respect, as the lifting numerical flux becomes central, with the difference that for BR2 the surface gradients are only penalized 
according to the local jump of the solution at the respective interface. Using BR2 rises the possibility to add an additional penalty to the gradients at wall 
boundaries. With the given coarse resolution the velocities at the wall would otherwise be unacceptable large. This is a result of the weak wall boundaries 
commonly used for DG schemes. We set $\eta_{BR2}=40$ in order to obtain small velocities at the wall.\\
A constant pressure source term $dp/dx$ is used to force the flow, which prescribes the Reynolds number of the flow exactly. The simulations are run until 
the bulk velocity statistically converged. Quantities plotted in 
\ref{fig:channel} are obtained by averaging 
over homogeneous directions.\\
The grid cells size is constant in x and z direction, in y direction a bell shape stretching is used with a ratio of $4$ from largest to smallest grid cell. In 
all directions $8$ cells are used.

\newpage
\bibliographystyle{acm}
\bibliography{References}

\begin{thebibliography}{10}

\bibitem{Bassi&Rebay:1997:B&F97}
{\sc Bassi, F., and Rebay, S.}
\newblock A high order accurate discontinuous finite element method for the
  numerical solution of the compressible {N}avier--{S}tokes equations.
\newblock {\em Journal of Computational Physics 131\/} (1997), 267--279.

\bibitem{bassi2}
{\sc Bassi, F., Rebay, S., Mariotti, G., Pedinotti, S., and Savini, M.}
\newblock A high-order accurate discontinuous finite element method for
  inviscid an viscous turbomachinery flows.
\newblock In {\em Proceedings of 2nd European Conference on Turbomachinery,
  Fluid and Thermodynamics\/} (Technologisch Instituut, Antwerpen, Belgium,
  1997), R.~Decuypere and G.~Dibelius, Eds., pp.~99--108.

\bibitem{les1}
{\sc Beck, A., Flad, D., Tonh\"auser, C., and Munz, C.-D.}
\newblock On the influence of polynomial de-aliasing on subgrid scale models.
\newblock {\em Flow, Turbulence and Combustion under revision\/} (2015).

\bibitem{les2}
{\sc Beck, A.~D., Bolemann, T., Flad, D., Frank, H., Gassner, G.~J.,
  Hindenlang, F., and Munz, C.-D.}
\newblock High-order discontinuous {Galerkin} spectral element methods for
  transitional and turbulent flow simulations.
\newblock {\em International Journal for Numerical Methods in Fluids 76}, 8
  (2014), 522--548.

\bibitem{carpenter2014}
{\sc Carpenter, M.~H., Fisher, T.~C., Nielsen, E.~J., and Frankel, S.~H.}
\newblock Entropy stable spectral collocation schemes for the navier--stokes
  equations: Discontinuous interfaces.
\newblock {\em SIAM Journal on Scientific Computing 36}, 5 (2014), B835--B867.

\bibitem{Carpenter-MH:2005kx}
{\sc Carpenter~MH, Kennedy~CA, B. H. V. S. V.~V.}
\newblock Fourth-order runge-kutta schemes for fluid mechanics applications.
\newblock {\em JOURNAL OF SCIENTIFIC COMPUTING 25}, 1 (October 2005), 157--194.

\bibitem{wiart2012}
{\sc {Carton de Wiart}, C., and Hillewaert, K.}
\newblock {DNS and ILES of transitional flows around a SD7003 using a high
  order discontinuous Galerkin method}.
\newblock In {\em Seventh International Conference on Computational Fluid
  Dynamics (ICCFD7)\/} (2012).

\bibitem{chandrashekar2012kinetic}
{\sc Chandrashekar, P.}
\newblock Kinetic energy preserving and entropy stable finite volume schemes
  for compressible euler and navier-stokes equations.
\newblock {\em arXiv preprint arXiv:1209.4994\/} (2012).

\bibitem{chollet1985}
{\sc Chollet, J.}
\newblock Two-point closure used for a sub-grid scale model in large eddy
  simulations.
\newblock In {\em Turbulent Shear Flows 4}. Springer, 1985, pp.~62--72.

\bibitem{desjardins2008high}
{\sc Desjardins, O., Blanquart, G., Balarac, G., and Pitsch, H.}
\newblock High order conservative finite difference scheme for variable density
  low mach number turbulent flows.
\newblock {\em Journal of Computational Physics 227}, 15 (2008), 7125--7159.

\bibitem{diosady2015}
{\sc Diosady, L., and Murman, S.}
\newblock Higher-order methods for compressible turbulent flows using entropy
  variables.
\newblock {\em AIAA Paper 294\/} (2015).

\bibitem{drikakis2009}
{\sc Drikakis, D., Hahn, M., Mosedale, A., and Thornber, B.}
\newblock Large eddy simulation using high-resolution and high-order methods.
\newblock {\em Philosophical Transactions of the Royal Society of London A:
  Mathematical, Physical and Engineering Sciences 367}, 1899 (2009),
  2985--2997.

\bibitem{fisher2013}
{\sc Fisher, T.~C., and Carpenter, M.~H.}
\newblock High-order entropy stable finite difference schemes for nonlinear
  conservation laws: Finite domains.
\newblock {\em Journal of Computational Physics 252\/} (2013), 518--557.

\bibitem{laf}
{\sc Flad, D., Beck, A., and Munz, C.-D.}
\newblock Simulation of underresolved turbulent flows by adaptive filtering
  using the high order discontinuous galerkin spectral element method.
\newblock {\em Journal of Computational Physics 313\/} (2016), 1--12.

\bibitem{gassner_skew_burgers}
{\sc Gassner, G.}
\newblock A skew-symmetric discontinuous {Galerkin} spectral element
  discretization and its relation to {SBP-SAT} finite difference methods.
\newblock {\em SIAM Journal on Scientific Computing 35}, 3 (2013),
  A1233--A1253.

\bibitem{tcfd2012}
{\sc Gassner, G., and Beck, A.}
\newblock On the accuracy of high-order discretizations for underresolved
  turbulence simulations.
\newblock {\em Theoretical and Computational Fluid Dynamics 27}, 3-4 (2013),
  221--237.

\bibitem{KoprivaGassner_Dispersion}
{\sc Gassner, G., and Kopriva, D.}
\newblock A comparison of the dispersion and dissipation errors of {Gauss} and
  {Gauss-Lobatto} discontinuous {Galerkin} spectral element methods.
\newblock {\em SIAM J. Scientific Computing 33}, 5 (2011), 2560--2579.

\bibitem{br1isstable}
{\sc Gassner, G., Winters, A., Hindenlang, F., and Kopriva, D.}
\newblock The br1 scheme is stable for the compressible navier-stokes
  equations.
\newblock {\em submitted to Journal of Scientific Computing\/}.

\bibitem{gassner2016split}
{\sc Gassner, G.~J., Winters, A.~R., and Kopriva, D.~A.}
\newblock Split form nodal discontinuous galerkin schemes with
  summation-by-parts property for the compressible euler equations.
\newblock {\em arXiv preprint arXiv:1604.06618\/} (2016).

\bibitem{germano91}
{\sc Germano, M., Piomelli, U., Moin, P., and Cabot, W.~H.}
\newblock A dynamic subgrid-scale eddy viscosity model.
\newblock {\em Physics of Fluids A: Fluid Dynamics (1989-1993) 3}, 7 (1991),
  1760--1765.

\bibitem{hickel2006}
{\sc Hickel, S., Adams, N.~A., and Domaradzki, J.~A.}
\newblock An adaptive local deconvolution method for implicit les.
\newblock {\em Journal of Computational Physics 213}, 1 (2006), 413--436.

\bibitem{VMShughes2000}
{\sc Hughes, T.~J., Mazzei, L., and Jansen, K.}
\newblock Large eddy simulation and the variational multiscale method.
\newblock {\em Computing and Visualization in Science 3}, 1-2 (2000), 47--59.

\bibitem{kennedy2008}
{\sc Kennedy, C., and Gruber, A.}
\newblock Reduced aliasing formulations of the convective terms within the
  {N}avier--{S}tokes equations for a compressible fluid.
\newblock {\em Journal of Computational Physics 227\/} (2008), 1676--1700.

\bibitem{Kirby2003}
{\sc Kirby, R., and Karniadakis, G.}
\newblock De-aliasing on non-uniform grids: algorithms and applications.
\newblock {\em Journal of Computational Physics 191\/} (2003), 249--264.

\bibitem{kraichnan1976eddy}
{\sc Kraichnan, R.~H.}
\newblock Eddy viscosity in two and three dimensions.
\newblock {\em Journal of the Atmospheric Sciences 33}, 8 (1976), 1521--1536.

\bibitem{kravchenko1997effect}
{\sc Kravchenko, A., and Moin, P.}
\newblock On the effect of numerical errors in large eddy simulations of
  turbulent flows.
\newblock {\em Journal of Computational Physics 131}, 2 (1997), 310--322.

\bibitem{lesieur1996}
{\sc Lesieur, M., and Metais, O.}
\newblock New trends in large-eddy simulations of turbulence.
\newblock {\em Annual review of fluid mechanics 28}, 1 (1996), 45--82.

\bibitem{lilly92}
{\sc Lilly, D.~K.}
\newblock A proposed modification of the germano subgrid-scale closure method.
\newblock {\em Physics of Fluids A: Fluid Dynamics (1989-1993) 4}, 3 (1992),
  633--635.

\bibitem{mittal1997suitability}
{\sc Mittal, R., and Moin, P.}
\newblock Suitability of upwind-biased finite difference schemes for large-eddy
  simulation of turbulent flows.
\newblock {\em AIAA journal 35}, 8 (1997), 1415--1417.

\bibitem{moser1999}
{\sc Moser, R.~D., Kim, J., and Mansour, N.~N.}
\newblock Direct numerical simulation of turbulent channel flow up to re= 590.
\newblock {\em Phys. Fluids 11}, 4 (1999), 943--945.

\bibitem{moura2016}
{\sc Moura, R., Mengaldo, G., Peir{\'o}, J., and Sherwin, S.}
\newblock On the eddy-resolving capability of high-order discontinuous galerkin
  approaches to implicit les/under-resolved dns of euler turbulence.
\newblock {\em Journal of Computational Physics\/} (2016).

\bibitem{moura2015}
{\sc Moura, R., Sherwin, S., and Peir{\'o}, J.}
\newblock Linear dispersion--diffusion analysis and its application to
  under-resolved turbulence simulations using discontinuous galerkin
  spectral/hp methods.
\newblock {\em Journal of Computational Physics 298\/} (2015), 695--710.

\bibitem{murman2016}
{\sc Murman, S.~M., Diosady, L.~T., Garai, A., and Ceze, M.}
\newblock A space-time discontinuous-galerkin approach for separated flows.
\newblock {\em Work 10\/} (2016), 5.

\bibitem{nagarajan2003robust}
{\sc Nagarajan, S., Lele, S.~K., and Ferziger, J.~H.}
\newblock A robust high-order compact method for large eddy simulation.
\newblock {\em Journal of Computational Physics 191}, 2 (2003), 392--419.

\bibitem{osswald2015l2roe}
{\sc O{\ss}wald, K., Siegmund, A., Birken, P., Hannemann, V., and Meister, A.}
\newblock L2roe: a low dissipation version of roe's approximate riemann solver
  for low mach numbers.
\newblock {\em International Journal for Numerical Methods in Fluids\/} (2015).

\bibitem{park2004discretization}
{\sc Park, N., Yoo, J.~Y., and Choi, H.}
\newblock Discretization errors in large eddy simulation: on the suitability of
  centered and upwind-biased compact difference schemes.
\newblock {\em Journal of Computational Physics 198}, 2 (2004), 580--616.

\bibitem{pirozzoli2011numerical}
{\sc Pirozzoli, S.}
\newblock Numerical methods for high-speed flows.
\newblock {\em Annual review of fluid mechanics 43\/} (2011), 163--194.

\bibitem{smago}
{\sc Smagorinsky, J.}
\newblock General circulation experiments with the primitive equations.
\newblock {\em Mon. Wea. Rev. 91\/} (1963), 99--164.

\bibitem{uranga2011}
{\sc Uranga, A., Persson, P.-O., Drela, M., and Peraire, J.}
\newblock Implicit large eddy simulation of transition to turbulence at low
  {R}eynolds numbers using a discontinuous {G}alerkin method.
\newblock {\em International Journal for Numerical Methods in Engineering 87},
  1-5 (2011), 232--261.

\bibitem{wiart2015}
{\sc Wiart, C.~C., Hillewaert, K., Bricteux, L., and Winckelmans, G.}
\newblock Implicit les of free and wall-bounded turbulent flows based on the
  discontinuous galerkin/symmetric interior penalty method.
\newblock {\em International Journal for Numerical Methods in Fluids 78}, 6
  (2015), 335--354.

\bibitem{yamazaki2002effects}
{\sc Yamazaki, Y., Ishihara, T., and Kaneda, Y.}
\newblock Effects of wavenumber truncation on high-resolution direct numerical
  simulation of turbulence.
\newblock {\em Journal of the Physical Society of Japan 71}, 3 (2002),
  777--781.

\end{thebibliography}

\end{document}